\documentclass[12pt,a4paper]{amsart}

\usepackage{amsmath,amsfonts,amsthm,amssymb,amscd,extarrows}

\usepackage[dvipdfmx]{graphicx}
\usepackage{here}   
\usepackage[stable]{footmisc}
\usepackage{cancel}
\usepackage{ascmac,fancybox}
\usepackage{cases}   
\usepackage{ulem}   
\usepackage[all]{xy}
\usepackage{url}
\usepackage{color}

\makeatletter
  
  \@addtoreset{equation}{section}
\makeatother
\newtheorem{Theorem}{Theorem}[section]
\newtheorem{Proposition}[Theorem]{Proposition}
\newtheorem{Lemma}[Theorem]{Lemma}

\newtheorem{Remark}[Theorem]{Remark}
\newtheorem{Example}[Theorem]{Example}
\newtheorem{Definition}[Theorem]{Definition}

\makeatletter
  
  \@addtoreset{equation}{section}
\makeatother

\def\RMN#1{\uppercase\expandafter{\romannumeral#1}}

\newcommand{\spec}{\mathop{\rm Spec}\nolimits}
\newcommand{\proj}{\mathop{\rm Proj}\nolimits}

\newcommand{\bfp}{\mathfrak{p}}

\newcommand{\bZ}{\mathbb{Z}}
\newcommand{\bR}{\mathbb{R}}
\newcommand{\bRo}{\mathbb{R}_{\ge 0}}

\newcommand{\bC}{\mathbb{C}}
\newcommand{\bQ}{\mathbb{Q}}
\newcommand{\bP}{\mathbb{P}}

\newcommand{\bNo}{\mathbb{N}_0}

\newcommand{\bma}{{\bf a}}
\newcommand{\bmb}{{\bf b}}
\newcommand{\bmx}{{\bf x}}
\newcommand{\bmp}{{\bf p}}

\newcommand{\cal}{\mathcal}

\newcommand{\Proof}{\noindent {\it Proof.} \ }
\newcommand{\seisei}[1]{\langle{#1}\rangle}

\begin{document}
\title[Equations of negative curves]{
Equations of negative curves of blow-ups of Ehrhart rings of rational convex polygons}
\author{Kazuhiko Kurano}
\date{}
\maketitle

\begin{abstract}
Finite generation of the symbolic Rees ring of a space monomial prime ideal of a $3$-dimensional weighted polynomial ring is a very interesting problem.
Negative curves play important roles in finite generation of these rings.
We are interested in the structure of the negative curve.
We shall prove that negative curves are rational in many cases.

We also see that the Cox ring of the blow-up of a toric variety at the point $(1,1,\ldots,1)$ coincides with the extended symbolic Rees ring of an ideal of a polynomial ring.
For example, Roberts' second counterexample to Cowsik's question (and Hilbert's 14th problem) coincides with the Cox ring of some normal projective variety (Remark~\ref{coxY}).
\end{abstract}

\section{Introduction}

Finite generation of symbolic Rees rings is one of very interesting problems in commutative ring theory.
It is deeply related to Hilbert's 14th problem or Kronecker's problem (Cowsik's question~\cite{Cowsik}).
It often happens that the Cox ring (or a multi-section ring) of an algebraic variety coincides with the extended symbolic Rees ring of an ideal of a ring.
Finite generation of these rings is also a very important problem in birational geometry.

The simplest non-trivial examples of symbolic Rees rings are 
\[
R_s(\bmp_{a,b,c}) = \bigoplus_{n \ge 0}{\bmp_{a,b,c}}^{(n)}t^n
 \subset k[x,y,z,t],
 \]
 where $\bmp_{a,b,c}$ is the kernel of the $k$-algebra homomorphism
\[
\phi_{a,b,c} : k[x,y,z] \longrightarrow k[\lambda]
\]
given by $\phi_{a,b,c}(x) = \lambda^a$, $\phi_{a,b,c}(y) = \lambda^b$, $\phi_{a,b,c}(z) = \lambda^c$
($k$ is a field and $a$, $b$, $c$ are pairwise coprime integers) and
${\bmp_{a,b,c}}^{(n)} = {\bmp_{a,b,c}}^{n}k[x,y,z]_{\bmp_{a,b,c}} \cap k[x,y,z]$ is the $n$th symbolic power.
In this case, the extended symbolic Rees ring $R_s(\bmp_{a,b,c})[t^{-1}]$ is the Cox ring of 
the blow-up $Y_{\Delta_{a,b,c}}$ of the weighted projective surface $X_{\Delta_{a,b,c}}=\proj(k[x,y,z])$ at the point $(1,1)$ in the torus.
Many commutative algebraists and algebraic geometers studied them and gave many results
(\cite{Hu}, \cite{C}, \cite{GNSnagoya}, \cite{GNW}, \cite{GK}, \cite{GAGK}, \cite{KM}, \cite{KN} etc.).
Finite generation of $R_s(\bmp_{a,b,c})$ depends on $a$, $b$, $c$ and the characteristic of $k$.
There are many examples of finitely generated $R_s(\bmp_{a,b,c})$.
Examples of infinitely generated $R_s(\bmp_{a,b,c})$ were first discovered by Goto-Nishida-Watanabe~\cite{GNW}.
In the case where the characteristic of $k$ is positive,
it is not known whether $R_s(\bmp_{a,b,c})$ is finitely generated for any $a$, $b$, $c$ or not.
We say that $C$ is a negative curve if $C$ is a curve in $Y_{\Delta_{a,b,c}}$ with $C^2 < 0$ such that
$C$ is not the exceptional curve $E$ of the blow-up $\pi : Y_{\Delta_{a,b,c}} \rightarrow X_{\Delta_{a,b,c}}$.
Since the Picard number of  $Y_{\Delta_{a,b,c}}$ is two, such a curve is unique if it exists. 
Cutkosky~\cite{C} proved that finite generation of the Cox ring of  $Y_{\Delta_{a,b,c}}$ is equivalent to the existence of curves $D_1$ and $D_2$ in $Y_{\Delta_{a,b,c}}$ such that $D_1 \cap D_2 = \emptyset$, $D_1\neq E$ and $D_2\neq E$ (the defining equations of  $\pi(D_1)$ and $\pi(D_2)$ satisfy Huneke's criterion~\cite{Hu} for finite generation).
If $R_s(\bmp_{a,b,c})$ is finitely generated with $\sqrt{abc} \not\in \bZ$,  
we may assume that either $D_1$ or $D_2$ coincides with the negative curve $C$.
If the negative curve does not exists, one can prove Nagata's conjecture (for $abc$ points) affirmatively  as in \cite{CK}.
Therefore the existence of the negative curve is a very important question.

The aim of this paper is to study the structure of the negative curve.
In particular, the author is interested in the problem whether the negative curve is rational or not.
If there exists a negative rational curve $C$, it is possible to estimate the degree of the curve $D$ 
such that $C$ and $D$ satisfies Huneke's criterion for finite generation in the same way as in \cite{KN}.

Assume that the negative curve $C$ in $Y_{\Delta_{a,b,c}}$ exists.
We shall study 
\[
\pi(C) \cap T
\]
in this paper, where $T$ is the torus in $X_{\Delta_{a,b,c}}$.
When $C.E = r$, the defining equation of $\pi(C) \cap T$ in $T = \spec(k[v^{\pm 1}, w^{\pm 1}])$
is an irreducible Laurent polynomial contained in $(v-1,w-1)^rk[v^{\pm 1}, w^{\pm 1}]$ and the Newton polygon has area less than $r^2/2$ (see Proposition~\ref{Prop3.2}).
We call such a Laurent polynomial an {\em $r$-nct} in this paper (Definition~\ref{Def3.1}).
Remark that $C$ and $\pi(C) \cap T$ is birational.
The author does not know any example that $C$ is not rational.

We shall prove some basic properties on $r$-ncts in Proposition~\ref{Fact3.3}.

It is proved that, for each $r \ge 0$, there exist essentially finitely many $r$-ncts.
When $r=1$, there exists essentially only one $1$-nct $\varphi_1 = vw-1$.
When $r=2$, there exists essentially only one $2$-nct $\varphi_2 = -v^2w -vw^2+3vw-1$.
When $r=3$, there exist essentially two $3$-ncts $\varphi_3$ and $\varphi'_3$ as in Example~\ref{123-nct}.
Let $P_{\varphi'_3}$ be the Newton polygon of $\varphi'_3$.
In the case where the characteristic of $k$ is not $2$, $P_{\varphi'_3}$ is a tetragon.
However, in the case of characteristic $2$, $P_{\varphi'_3}$ is a smaller triangle than this tetragon.
This is a reason why ${\bmp_{9,10,13}}^{(3)}$ contains a negative curve in the case of characteristic $2$. 
In other characteristic, ${\bmp_{9,10,13}}^{(3)}$ does not contain a negative curve (see Example~\ref{ExampleOfNct}).

The following is the main theorem of this paper.

\begin{Theorem}\label{Thm3.6}
\begin{rm}
Let $k$ be a field and $\varphi$ be an $r$-nct over $k$, where $r \ge 2$.
Let $P_\varphi$ be the Newton polygon of $\varphi$.
Consider the prime ideal $\bmp_\varphi$ of the Ehrhart ring 
\[
E(P_\varphi,\lambda) = \bigoplus_{d \ge 0}
\left( \bigoplus_{(\alpha,\beta) \in d P_\varphi \cap \bZ^2} kv^\alpha w^\beta \right) \lambda^d
\subset k[v^{\pm 1}, w^{\pm 1}][\lambda]
\]
satisfying $\bmp_\varphi = E(P_\varphi,\lambda) \cap (v-1,w-1)k[v^{\pm 1}, w^{\pm 1}][\lambda]$.
Let $Y_{\Delta_\varphi}$ be the blow-up of $X_{\Delta_\varphi} = \proj(E(P_\varphi, \lambda))$ at the point $\bmp_\varphi$.
Let $C_\varphi$ be the proper transform of the curve $V_+(\varphi \lambda)$ in $X_{\Delta_\varphi}$.
Let $I_\varphi$ (resp.\ $B_\varphi$) be the number of the interior lattice points 
(resp.\ the boundary lattice points) of $P_\varphi$.

Consider the following conditions:
\begin{enumerate}
\item[(1)]
$-K_{Y_{\Delta_\varphi}}$ is nef and big.
\item[(2)]
$-K_{Y_{\Delta_\varphi}}$ is nef.
\item[(3)]
$(-K_{Y_{\Delta_\varphi}})^2>0$.
\item[(4)]
Let $\bma_1$, $\bma_2$, \ldots, $\bma_n$ be the first lattice points of the $1$-dimensional cones in the fan corresponding to the toric variety $X_{\Delta_\varphi}$.
We put
\[
P_{-K_{X_{\Delta_\varphi}}} =
\{ \bmx \in \bR^2 \mid \seisei{\bmx, \bma_i} \ge -1 \ (i = 1, 2, \ldots, n) \} .
\]
Then $|P_{-K_{X_{\Delta_\varphi}}}| > 1/2$ is satisfied.
\item[(5)]
$-K_{Y_{\Delta_\varphi}}$ is big.
\item[(6)]
The Cox ring ${\rm Cox}(Y_{\Delta_\varphi})$ is Noetherian.
\footnote{
As we shall see in Remark~\ref{coxY},
${\rm Cox}(Y_{\Delta_\varphi})$ is a extended symbolic Rees ring of an ideal $I$ over
the polynomial ring ${\rm Cox}(X_{\Delta_\varphi})$.
If ${\rm Cl}(X_{\Delta_\varphi})$ is torsion-free, then $I$ is a prime ideal.
In the case where the characteristic of $k$ is $0$, it holds $\sqrt{I} = I$.}
\item[(7)]
$B_\varphi \ge r$.
\item[(8)]
$H^0(Y_{\Delta_\varphi}, {\cal O}_{Y_{\Delta_\varphi}}(K_{Y_{\Delta_\varphi}} + n C_\varphi)) = 0$ for $n > 0$. 
\item[(9)]
$I_\varphi = \frac{r(r-1)}{2}$.
\item[(10)]
The extended symbolic Rees ring $R'_s(\bmp_\varphi)$ is Noetherian.
\item[(11)]
$C_\varphi \simeq \bP_k^1$ is satisfied.
\end{enumerate}

Then we have the following:
\begin{itemize}
\item[(a)]
\[
\begin{array}{ccccccccccc}
(1) & \Longrightarrow & (2) & \multicolumn{3}{c}{\xLongrightarrow[]{\hspace{5.5em}}}  &(7) & & & & \\
\Downarrow & & & & & & \Downarrow& &  & & \\
(3) & \Longrightarrow & (4) & \Longrightarrow & (5) & \Longrightarrow & (8) & \Longleftrightarrow & (9) & \Longleftarrow & (11) \\
& & & & \Downarrow & & \Downarrow& & & & \\
& & & & (6) & \Longrightarrow & (10) & & & & 
\end{array}
\]
is satisfied.\footnote{
We do not have to restrict ourselves to the polygon $P_\varphi$ for 
$(1)\Rightarrow (3)\Rightarrow  (4)\Rightarrow (5)\Rightarrow  (6)\Rightarrow  (10)$.
They hold for any integral convex polygon.}
\item[(b)]
If $k$ is algebraically closed, then (9) is equivalent to (11).
\item[(c)]
If the characteristic of $k$ is positive, then (10) is always satisfied.
\end{itemize}
\end{rm}
\end{Theorem}

We shall prove this theorem in section~\ref{sect4}.
Using this theorem, we shall see that the negative curve $C$ in $Y_{\Delta_{a,b,c}}$ is rational in many cases in Remark~\ref{rational}.
For instance, if $r \le 4$, then $C$ is rational.
If $r \le 5$ and ${\rm ch}(k) = 0$, $C$ is rational.

\section{Preliminaries}\label{junbi}

All rings in this paper are commutative with unity.

For a Noetherian ring $A$ and a prime ideal $Q$ of $A$, 
$Q^{(n)} = Q^nA_Q \cap A$ is called the $n$th symbolic power of $Q$.
When $n < 0$, we think $Q^{(n)} = A$.
We put
\begin{eqnarray*}
R_s(Q) & = & \oplus_{n \ge 0} Q^{(n)} t^n \subset A[t] \\
R'_s(Q) & = & \oplus_{n \in \bZ} Q^{(n)} t^n \subset A[t^{\pm 1}]
\end{eqnarray*}
and call them the symbolic Rees ring of $Q$ and
the extended symbolic Rees ring of $Q$, respectively.
We know that $R_s(Q)$ is Noetherian if and only if so is $R'_s(Q)$.

\begin{Definition}\label{Def1}
\begin{rm}
Let $a$, $b$, $c$ be pairwise coprime integers.
Let $k$ be a field and $S_{a,b,c} = k[x,y,z]$ be a graded polynomial ring with $\deg(x) = a$, $\deg(y) = b$, $\deg(z) = c$.
Let $\lambda$ be a variable and consider the $k$-algebra homomorphism
\[
\phi_{a,b,c} : S_{a,b,c} \longrightarrow k[\lambda]
\]
given by $\phi_{a,b,c}(x) = \lambda^a$, $\phi_{a,b,c}(y) = \lambda^b$, $\phi_{a,b,c}(z) = \lambda^c$.
Let $\bmp_{a,b,c}$ be the kernel of $\phi_{a,b,c}$.
\end{rm}
\end{Definition}

\begin{Remark}\label{Rem2}
\begin{rm}
Let $Y_{\Delta_{a,b,c}}$ be the blow-up of $X_{\Delta_{a,b,c}} = \proj S_{a,b,c}$ at $V_+(\bmp_{a,b,c})$.
(Here $X_{\Delta_{a,b,c}}$ is a toric variety with a fan $\Delta_{a,b,c}$.
We shall define $\Delta_{a,b,c}$ in Remark~\ref{MethodOfToric}.)
Let $E$ be the exceptional divisor and $H$ be the pull back of ${\cal O}_{X_{\Delta_{a,b,c}}}(1)$.
Then $\{ E, H \}$ is a generating set of ${\rm Cl}(Y_{\Delta_{a,b,c}}) \simeq \bZ^2$.
By Cutkosky~\cite{C}, we have an identification
\[
H^0(Y_{\Delta_{a,b,c}}, {\cal O}_{Y_{\Delta_{a,b,c}}}(dH - rE)) = [{\bmp_{a,b,c}}^{(r)}]_d
\]
and an isomorphism of rings as follows:
\begin{equation}\label{doukei}
{\rm Cox}(Y_{\Delta_{a,b,c}}) = \bigoplus_{d, r \in \bZ}H^0(Y_{\Delta_{a,b,c}}, {\cal O}_{Y_{\Delta_{a,b,c}}}(dH - rE))
= \bigoplus_{r \in \bZ} {\bmp_{a,b,c}}^{(r)} t^r = R'_s(\bmp_{a,b,c}) 
\end{equation}
\end{rm}
\end{Remark}

\begin{Remark}\label{MethodOfToric}
\begin{rm}
Consider $S_{a,b,c}$, $\bmp_{a,b,c}$, $X_{\Delta_{a,b,c}}$, $Y_{\Delta_{a,b,c}}$ as in 
Definition~\ref{Def1} and Remark~\ref{Rem2}.

By Herzog~\cite{Her}, we have
\[
\bmp_{a,b,c} = I_2\left(
\begin{array}{ccc}
x^{s_2} & y^{t_3} & z^{u_1} \\ y^{t_1} & z^{u_2} & x^{s_3}
\end{array}
\right)
= (x^{s} - y^{t_1}z^{u_1}, y^{t} - x^{s_2}z^{u_2}, z^{u} - x^{s_3}y^{t_3}) ,
\]
where $s = s_2+s_3$, $t = t_1 + t_3$, $u = u_1+u_2$, and
$I_2( \ )$ is the ideal generated by $2 \times 2$-minors of the given matrix.
Here we put $v = x^{s_2}z^{u_2}/y^t$, $w = x^{s_3}y^{t_3}/z^u$.
Then $S_{a,b,c}[x^{-1}, y^{-1}, z^{-1}]$ is a ${\Bbb Z}$-graded ring such that
\[
S_{a,b,c}[x^{-1}, y^{-1}, z^{-1}]_0=k[v^{\pm 1}, w^{\pm 1}] 
\]
(cf.\  the proof of Lemma~3.6 in \cite{KN}).

Take integers $i_0$, $j_0$ satisfying $i_0a + j_0b = 1$.
Putting $\lambda = x^{i_0}y^{j_0}$, we obtain
\[
S_{a,b,c} \subset S_{a,b,c}[x^{-1}, y^{-1}, z^{-1}]= \left( S_{a,b,c}[x^{-1}, y^{-1}, z^{-1}]_0 \right) [\lambda^{\pm 1}]
= k[v^{\pm 1}, w^{\pm 1}, \lambda^{\pm 1}] .
\]
Here suppose
\[
v^\alpha w^\beta \lambda^d  =  x^iy^jz^k \in S_{a,b,c} .
\]
Then we have $d = \deg(x^iy^jz^k) = ia + jb + kc \ge 0$ and
\begin{align*}
& i = s_2\alpha + s_3 \beta + i_0d \ge 0 \\
& j = -t\alpha + t_3 \beta + j_0d \ge 0 \\
& k = u_2\alpha - u \beta  \ge 0  .
\end{align*}
Consider the rational triangle
\[
P_{a,b,c} = \left\{ (\alpha,\beta) \in \bR^2 \ \left| \ 
\begin{array}{l}
s_2\alpha + s_3 \beta + i_0 \ge 0 \\
-t\alpha + t_3 \beta + j_0 \ge 0 \\
u_2\alpha - u \beta  \ge 0
\end{array}
\right\} \right. .
\]
When $s_3$, $t_3$, $u$ are positive integers\footnote{
If $\bmp_{a,b,c}$ is not a complete intersection,
then all of $s_2$, $s_3$, $t_1$, $t_3$, $u_1$, $u_2$ are positive by Herzog~\cite{Her}.
In the case where $\bmp_{a,b,c}$ is a complete intersection,
we can choose positive integers $s_3$, $t_3$, $u$ after a suitable permutation of $a$, $b$, $c$.
},
$P_{a,b,c}$ is the following triangle:
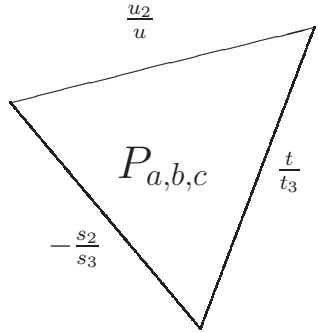
\[
{
\setlength\unitlength{1truecm}
  \begin{picture}(6,5)(-1,-3)
  \put(1.4,-1){\mbox{\Large $P_{a,b,c}$}}
  \put(0,0){\line(4,1){4}}
\qbezier (0,0) (1.25,-1.5) (2.5,-3)
\qbezier (4,1) (3.25,-1) (2.5,-3)
  \put(1.5,1){$\frac{u_2}{u}$}
  \put(3.5,-1){$\frac{t}{t_3}$}
  \put(0.5,-2){$-\frac{s_2}{s_3}$}
  \end{picture}
}
\]
Thus we have the identification
\begin{equation}\label{iden}
(S_{a,b,c})_d = \left( \bigoplus_{(\alpha, \beta) \in d P_{a,b,c}  \cap \bZ^2} k v^\alpha w^\beta \right) \lambda^d .
\end{equation}
Consider the Ehrhart ring of $P_{a,b,c}$ defined by 
\[
E(P_{a,b,c}, \lambda) = \bigoplus_{d \ge 0}\left( \bigoplus_{(\alpha, \beta) \in dP_{a,b,c} \cap \bZ^2} k v^\alpha w^\beta \right) \lambda^d \subset k[v^{\pm 1}, w^{\pm 1}, \lambda^{\pm 1}] .
\]
Thus $E(P_{a,b,c}, \lambda)$ is isomorphic to $S_{a,b,c}$ as a graded $k$-algebra.
Let $\Delta_{a,b,c}$ be the complete fan in $\bR^2$ with one dimensional cones
\[
\bRo (s_2,s_3), \ \ \bRo (-t,t_3), \ \ \bRo (u_2,-u) .
\]
Then $\proj(E(P_{a,b,c}, \lambda))$ is the toric variety $X_{\Delta_{a,b,c}}$ with the fan
$\Delta_{a,b,c}$.

Furthermore we have 
\[
\bmp_{a,b,c} = E(P_{a,b,c}, \lambda) \cap (v-1, w-1)k[v^{\pm 1}, w^{\pm 1}, \lambda^{\pm 1}] .
\]
It is easy to see 
\begin{equation}\label{nth}
{\bmp_{a,b,c}}^{(r)} = E(P_{a,b,c}, \lambda) \cap (v-1, w-1)^rk[v^{\pm 1}, w^{\pm 1}, \lambda^{\pm 1}]
\end{equation}
for any $r \ge 1$.
Therefore we have
\[
H^0(Y_{\Delta_{a,b,c}}, {\cal O}_{Y_{\Delta_{a,b,c}}}(dH - rE)) = 
\left[ \bfp^r \cap 
\left( \bigoplus_{(\alpha, \beta) \in d P_{a,b,c}  \cap \bZ^2} k v^\alpha w^\beta \right)
\right] \lambda^d 
 = [{\bmp_{a,b,c}}^{(r)}]_d
,
\]
where $\bfp =  (v-1, w-1)k[v^{\pm 1}, w^{\pm 1}]$.
\end{rm}
\end{Remark}

\begin{Remark}\label{Area}
\begin{rm}
With notation as in Remark~\ref{MethodOfToric}, we shall calculate the area of $P_{a,b,c}$ here.

Suppose $z^u = v^{\alpha_0}w^{\beta_0}\lambda^{cu}$.
Then the bottom vertex of $cuP_{a,b,c}$ is $(\alpha_0,\beta_0)$.
Since
\[
x^{s_3}y^{t_3} = wz^u = v^{\alpha_0}w^{\beta_0+1}\lambda^{cu} ,
\]
the point $(\alpha_0,\beta_0+1)$ is on the upper edge of $cuP_{a,b,c}$.
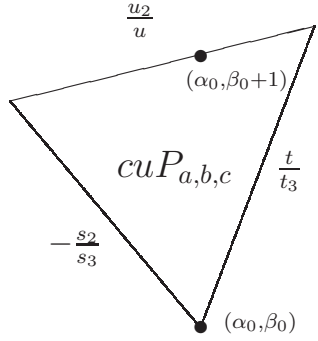
\[
{
\setlength\unitlength{1truecm}
  \begin{picture}(6,5)(-1,-3)
  \put(2.8,-3){${\scriptstyle (\alpha_0,\beta_0)}$}
  \put(2.3,0.2){${\scriptstyle (\alpha_0,\beta_0+1)}$}
  \put(1.4,-1){\mbox{\large $cuP_{a,b,c}$}}
  \put(0,0){\line(4,1){4}}
\qbezier (0,0) (1.25,-1.5) (2.5,-3)
\qbezier (4,1) (3.25,-1) (2.5,-3)
  \put(1.5,1){$\frac{u_2}{u}$}
  \put(3.5,-1){$\frac{t}{t_3}$}
  \put(0.5,-2){$-\frac{s_2}{s_3}$}
  \put(2.4,-3.1){$\bullet$}
  \put(2.4,0.5){$\bullet$}
  \end{picture}
}
\]
Then the width of $cuP_{a,b,c}$ is 
\[
\frac{1}{\frac{t}{t_3} - \frac{u_2}{u}} + \frac{1}{\frac{s_2}{s_3} + \frac{u_2}{u}}
= \frac{t_3 u}{a} + \frac{s_3u}{b} = \frac{cu^2}{ab} .
\]
Here recall that 
\begin{align*}
a & = e((\lambda^a), k[\lambda]) = e((x), S_{a,b,c}/\bmp_{a,b,c}) = 
\ell_{S_{a,b,c}}(S_{a,b,c}/(x)+\bmp_{a,b,c}) \\
& = \ell_{S_{a,b,c}}(S_{a,b,c}/(x, y^{t_1}z^{u_1}, y^t, z^u))
= tu - t_3u_2 ,
\end{align*}
$b = su - s_3u_1 = s_2u + s_3u_2$ and $cu = s_3a+t_3b$,
where $e( \ )$ is the multiplicity and $\ell( \ )$ is the length.
Therefore the area of $cuP_{a,b,c}$ is $\frac{cu^2}{2 ab}$.
Thus we know that the area of $P_{a,b,c}$ is $\frac{1}{2 abc}$.
\end{rm}
\end{Remark}

\begin{Definition}\label{negative curve}
\begin{rm}
If an irreducible polynomial $f \in [{\bmp_{a,b,c}}^{(r')}]_{d'}$ satisfies $d'/r' < \sqrt{abc}$,
we say that $f \in [{\bmp_{a,b,c}}^{(r')}]_{d'}$ is a negative curve.
\end{rm}
\end{Definition}

If a negative curve $f \in [{\bmp_{a,b,c}}^{(r')}]_{d'}$ exists,
then both $r'$ and $d'$ are uniquely determined, and $f$ is also unique
up to a constant factor.
We denote the proper transform of $V_+(f)$ by $C$ ($\subset Y_{\Delta_{a,b,c}}$).
Then $C$ satisfies $C^2 = \frac{{d'}^2}{abc} - {r'}^2 < 0$.

\begin{Lemma}\label{PrimeElementOfE}
Let $k$ be a field and $P$ be a rational convex polygon (the convex hull of a finite subset of $\bQ^2$) in $\bR^2$.
Let $A = k[v^{\pm 1}, w^{\pm 1}]$ be a Laurent polynomial ring with two variables $v$, $w$.
Consider the Ehrhart ring
\[
E(P, \lambda) = \bigoplus_{d \ge 0} \left( \bigoplus_{(\alpha,\beta) \in dP \cap \bZ^2}
kv^\alpha w^\beta \right) \lambda^d \subset A[\lambda] .
\]
Let $d_1$ be a positive integer and 
put $\varphi = \sum_{(\alpha,\beta) \in d_1P \cap \bZ^2}
c_{(\alpha,\beta)}v^\alpha w^\beta$ where $c_{(\alpha,\beta)} \in k$ and
\begin{equation}\label{Nvarphi}
N_\varphi = \{ (\alpha,\beta) \in \bZ^2 \mid c_{(\alpha,\beta)} \neq 0 \} .
\end{equation}

Assume that $N_\varphi$ contains at least two elements.
Then the following two conditions are equivalent:
\begin{enumerate}
\item
$\varphi \lambda^{d_1}$ is a prime element of $E(P, \lambda)$.
\item
$\varphi$ is irreducible in $A$ and each edge of $P$ contains an element in $N_\varphi$.
\end{enumerate}
\end{Lemma}

\Proof
Before proving this lemma, we shall give some remarks here.
Suppose that a lattice point $(\alpha', \beta')$ is in the interior of $d'P$ for some $d'>0$.
Then we have 
\[
E(P, \lambda)[(v^{\alpha'} w^{\beta'} \lambda^{d'})^{-1}] = A[\lambda^{\pm 1}]
\]
and
\begin{equation}\label{domain}
\left( E(P, \lambda)/(\varphi \lambda^{d_1}) \right)[(v^{\alpha'} w^{\beta'} \lambda^{d'})^{-1}]  
= \left( A/ \varphi A \right) [\lambda^{\pm 1}] .
\end{equation}
Here remark that both sides are not $0$ since $N_\varphi$ contains at least two elements.
Let $\{ E_1, \ldots, E_s \}$ be the set of edges of $P$.
Putting
\[
\bmp_i = \bigoplus_{d > 0} \left( \bigoplus_{(\alpha,\beta) \in d(P \setminus E_i) \cap \bZ^2}
kv^\alpha w^\beta \right) \lambda^d
\subset E(P, \lambda) ,
\]
$\bmp_i$ is a prime ideal of $E(P, \lambda)$ of height $1$ for $i = 1, 2, \ldots, s$.
It is easy to see 
\begin{equation}\label{radical}
\sqrt{(v^{\alpha'} w^{\beta'} \lambda^{d'})} = \bmp_1 \cap \cdots  \cap \bmp_s .
\end{equation}

First assume (1).
Since $\varphi \lambda^{d_1}$ is a prime element, $\varphi$ is irreducible in $A$ by (\ref{domain}).
If some $E_i$ does not meet $N_\varphi$, $\varphi \lambda^{d'}$ is contained in $\bmp_i$.
Since $\varphi \lambda^{d_1}$ is a prime element, we obtain $\bmp_i = (\varphi \lambda^{d_1})$.
It is impossible since $N_\varphi$ contains at least two points.

Next assume (2).
Since $\varphi$ is irreducible in $A$, the right-hand side of (\ref{domain}) is an integral domain.
In order to prove (1), it is enough to show that $\varphi \lambda^{d_1}$, $v^{\alpha'} w^{\beta'} \lambda^{d'}$
 is an $E(P, \lambda)$-regular sequence.
Since $E(P, \lambda)$ is a normal domain, we know
\[
{\rm Ass}_{E(P, \lambda)}\left( E(P, \lambda)/(v^{\alpha'} w^{\beta'} \lambda^{d'}) \right) = \{ \bmp_1, \ldots, \bmp_s \} 
\]
by (\ref{radical}).
Since none of $\bmp_i$'s contains $\varphi \lambda^{d_1}$ by (2),
$v^{\alpha'} w^{\beta'} \lambda^{d'}$, $\varphi \lambda^{d_1}$
is an $E(P, \lambda)$-regular sequence.
Since $v^{\alpha'} w^{\beta'} \lambda^{d'}$ and $\varphi \lambda^{d_1}$ are homogeneous elements,
$\varphi \lambda^{d_1}$, $v^{\alpha'} w^{\beta'} \lambda^{d'}$
 is an $E(P, \lambda)$-regular sequence.
\qed

\begin{Remark}\label{coxY}
\begin{rm}
As in Remark~\ref{Rem2}, the extended symbolic Rees ring $R'_s(\bmp_{a,b,c})$ 
is identified with the Cox ring of a blow-up of a toric surface.
Here we generalize this identification.

Let $X_\Delta$ be a $d$-dimensional toric variety with a fan $\Delta$.
Let 
\[
\{ \bRo \bma_1, \bRo \bma_2, \ldots, \bRo \bma_n \}
\]
be the set of the $1$-dimensional cones in $\Delta$.
We assume that each $\bma_i$ is in $\bZ^d$ such that the greatest common measure
of the components of $\bma_i$ is $1$, and $\sum_i\bR \bma_i = \bR^d$.
Here remark that $\bma_i \neq \bma_j$ if $i \neq j$.
Then we have the exact sequence 
\[
0 \longleftarrow {\rm Cl}(X_\Delta) \longleftarrow \bZ^n \stackrel{
\left(
\begin{array}{c}
\bma_1 \\ \bma_2 \\ \vdots \\ \bma_n
\end{array}
\right)
}{\longleftarrow}
\bZ^d \longleftarrow 0 ,
\]
where ${\rm Cl}(X_\Delta)$ is the divisor class group of $X_\Delta$.

The morphism of monoids
\begin{equation}\label{grading}
{\rm Cl}(X_\Delta) \longleftarrow \bZ^n \supset (\bNo)^n
\end{equation}
induces the morphism of semi-group rings
\[
k[{\rm Cl}(X_\Delta)] \longleftarrow A:= k[x_1, x_2, \ldots, x_n] = k[(\bNo)^n] .
\]
Let $I$ be the kernel of this ring homomorphism.
By (\ref{grading}), $A$ has a structure of a ${\rm Cl}(X_\Delta)$-graded ring.
Then we have an isomorphism $A \simeq {\rm Cox}(X_\Delta)$ as a ${\rm Cl}(X_\Delta)$-graded ring~\cite{Cox}.

Let $\{ \bmp_1, \bmp_2, \ldots, \bmp_\ell \}$ be the set of minimal prime ideals of $I$ and put
\[
I^{(r)} = I^rA_{\bmp_1} \cap \cdots \cap I^rA_{\bmp_\ell} \cap A .
\]
In the case where ${\rm Cl}(X_\Delta)$ is torsion-free, $I$ is a prime ideal of $A$
and $I^{(r)}$ is the $r$th symbolic power of $I$.
If the characteristic of $k$ is $0$, then we have $I = \sqrt{I}$ and
$I^{(r)} =\bmp_1^{(r)} \cap \bmp_2^{(r)} \cap \cdots \cap \bmp_\ell^{(r)}$.

Let $Y_\Delta$ be the blow-up of $X_\Delta$ at $(1,1,\dots,1)$ in the torus $(k^\times)^n$.
Then we have an isomorphism 
\begin{equation}\label{generalization}
{\rm Cox}(Y_\Delta) \simeq R'_s(I)
\end{equation}
as a ${\rm Cl}(X_\Delta) \times \bZ$-graded ring.

In the case of $X_{\Delta_{a,b,c}}$ and $Y_{\Delta_{a,b,c}}$ in Definition~\ref{Def1},  Remark~\ref{Rem2}
and Remark~\ref{MethodOfToric},
we put $\bma_1 = (s_2,s_3)$, $\bma_2 = (-t,t_3)$, $\bma_3= (u_2,-u)$.
Then the morphism (\ref{grading}) of monoids coincides with
\[
\bZ \stackrel{
(a \ b \ c)}{\longleftarrow} \bZ^3 \supset (\bNo)^3 .
\]
Let $I$ be the kernel of the $k$-algebra homomorphism $k[x_1,x_2,x_3] \rightarrow k[y_1^{\pm 1}]$ 
($x_1 \mapsto y_1^{a}$, $x_2 \mapsto y_1^b$, $x_3 \mapsto y_1^c$).
Then $I = \bmp_{a,b,c}$ and ${\rm Cox}(Y_{\Delta_{a,b,c}})\simeq R'_s(\bmp_{a,b,c})$ by (\ref{generalization}).
It is just the isomorphism in Remark~\ref{Rem2} given by Cutkosky.

Suppose $\bma_1 = (2,-1)$, $\bma_1 = (-2,-1)$, $\bma_1 = (0,1)$.
Then the morphism (\ref{grading}) of  monoids is
\[
\bZ \oplus \bZ/(2) 
\longleftarrow \bZ^3 \supset (\bNo)^3
\]
such that $(1,0,0) \mapsto (1,\overline{1})$, $(0,1,0) \mapsto (1,\overline{0})$, $(0,0,1) \mapsto (2,\overline{1})$, respectively.
Letting $I$ be the kernel of the $k$-algebra homomorphism $k[x_1,x_2,x_3] \rightarrow \frac{k[y_1^{\pm 1}, y_2]}{(y_2^2-1)}$ 
($x_1 \mapsto y_1y_2$, $x_2 \mapsto y_1$, $x_3 \mapsto y_1^2y_2$),
we have ${\rm Cox}(Y_\Delta)\simeq R'_s(I)$.
Consider the following $k$-algebra homomorphisms $f_1 : k[x_1,x_2,x_3] \rightarrow k[y_1]$ 
($x_1 \mapsto y_1$, $x_2 \mapsto y_1$, $x_3 \mapsto y_1^2$) and
$f_2 : k[x_1,x_2,x_3] \rightarrow k[y_1]$ 
($x_1 \mapsto -y_1$, $x_2 \mapsto y_1$, $x_3 \mapsto -y_1^2$).
If the characteristic of $k$ is not $2$, then we have $I = {\rm Ker}(f_1) \cap {\rm Ker}(f_2)$.

Suppose $t \ge 2$ and 
\[
\left(
\begin{array}{c}
\bma_1 \\ \bma_2 \\ \bma_3 \\ \bma_4 \\ \bma_5 \\ \bma_6 \\ \bma_7
\end{array}
\right)
= 
\left(
\begin{array}{ccc}
-1 & t & t \\
t & -1 & t \\
t & t & -1 \\
1 & 0 & 0 \\
0 & 1 & 0 \\
0 & 0 & 1 \\
-1 & -1 & -1
\end{array}
\right) .
\]
Then the morphism (\ref{grading}) of  monoids is
\[
\bZ^4 \stackrel{
\left(
\begin{array}{ccccccc}
1 & 0 & 0 & t+1 & 0 & 0 & t \\
0 & 1 & 0 & 0 & t+1 & 0 & t \\
0 & 0 & 1 & 0 & 0 & t+1 & t \\
0 & 0 & 0 & 1 & 1 & 1 & 1
\end{array}
\right)}{
\xlongleftarrow[]{\hspace{16em}}} \bZ^7 \supset (\bNo)^7 .
\]
Let $I$ be the kernel of the $k$-algebra homomorphism $k[x_1,x_2,x_3,S,T,U,V] \rightarrow k[x_1,x_2,x_3,W]$ 
($S \mapsto x_1^{t+1}W$, $T \mapsto x_2^{t+1}W$, $U \mapsto x_3^{t+1}W$,
$V \mapsto (x_1x_2x_3)^{t}W$).
Then $I$ is a prime ideal of height $3$ and ${\rm Cox}(Y_\Delta) \simeq R'_s(I)$.
Roberts~\cite{Ro} proved that $R'_s(I)$ is not Noetherian and
some pure subring of $R'_s(I)$ gives a counterexample to Hilbert's fourteenth problem
if the characteristic of $k$ is $0$.
On the other hand, $R'_s(I)$ is Noetherian if  the characteristic of $k$ is positive~\cite{K8}.

Let $E$ be the Ehrhart ring of the following tetragon
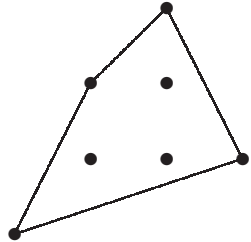
\[
{
\setlength\unitlength{1truecm}
  \begin{picture}(6,3.5)(-2,-2)
\qbezier (0,-2) (1.5,-1.5) (3,-1)
\qbezier (0,-2) (0.5,-1) (1,0)
\qbezier (1,0) (1.5,0.5) (2,1)
\qbezier (2,1) (2.5,0) (3,-1)
\put(-0.1,-2.1){$\bullet$}
\put(0.9,-0.1){$\bullet$}
\put(0.9,-1.1){$\bullet$}
\put(1.9,-0.1){$\bullet$}
\put(1.9,-1.1){$\bullet$}
\put(1.9,0.9){$\bullet$}
\put(2.9,-1.1){$\bullet$}
  \end{picture}
}
\]
and put $X_\Delta = \proj(E)$.
In this case, the $1$-dimensional cones are
\[
\left(
\begin{array}{c}
\bma_1 \\ \bma_2 \\ \bma_3 \\ \bma_4
\end{array}
\right)
= 
\left(
\begin{array}{ccc}
2 & -1 \\
1 & -1 \\
-2 & -1 \\
-1 & 3 
\end{array}
\right) 
\]
and   the morphism (\ref{grading}) of  monoids is
\[
\bZ^2 \stackrel{
\left(
\begin{array}{cccc}
1 & 1 & 1 & 1 \\ 3 & -4 & 1 & 0
\end{array}
\right)}{
\longleftarrow} \bZ^4 \supset (\bNo)^4 .
\]
Let $I$ be the kernel of the $k$-algebra homomorphism
 $k[x_1,x_2,x_3,x_4] \rightarrow k[y_1^{\pm 1},y_2^{\pm 1}]$ 
($x_1 \mapsto y_1y_2^3$, $x_2 \mapsto y_1y_2^{-4}$, $x_3 \mapsto y_1y_2$, $x_4 \mapsto y_1$).
Then $I$ is a prime ideal called a toric ideal.
By Example~\ref{ExampleOfGGK},
we know that $R'_s(I)$ is a Noetherian ring.
Let $\bmp$ be a prime ideal of the Ehrhart ring $E \subset k[v^{\pm 1}, w^{\pm 1}][\lambda]$ satisfying $\bmp = E \cap (v-1,w-1)k[v^{\pm 1}, w^{\pm 1}][\lambda]$.
Then the symbolic Rees ring $R'_s(\bmp)$ is a pure subring of $R'_s(I)$.
Therefore we know that $R'_s(\bmp)$ is also a Noetherian ring.
\end{rm}
\end{Remark}

\section{Definition of $r$-nct and basic properties}\label{nctdefbasic}

We define an $r$-nct and give basic properties of it in this section.

\begin{Definition}\label{Def3.1}
\begin{rm}
Let $k$ be a field and $A = k[v^{\pm 1}, w^{\pm 1}]$ be a Laurent polynomial ring with two variables $v$, $w$.
Let $r$ be a positive integer.
An element $\varphi$ of $A$ is called an {\it equation of a negative curve with multiplicity $r$ over $k$} (or simply an {\it $r$-nct over $k$})
if the following two conditions are satisfied:
\begin{enumerate}
\item[(1)]
$\varphi$ is irreducible in $A$ and $\varphi \in (v-1,w-1)^rA$.
\item[(2)]
Let $N_\varphi$ be the set defined as in (\ref{Nvarphi}).
Let $P_\varphi$ be the convex hull of $N_\varphi$.
Then $2|P_\varphi| < r^2$ is satisfied, where $|P_\varphi|$ is the area of the polygon $P_\varphi$.
\end{enumerate}
\end{rm}
\end{Definition}

\begin{Example}\label{123-nct}
\begin{rm}
Polynomials $v-1$, $w-1$, $\varphi_1=vw-1$ are $1$-ncts over any field $k$.
The area of $P_{\varphi_1}$ is $0$.

The polynomial
\[
\varphi_2=-\varphi_1(v-1) - v(w-1)^2 = -v^2w - vw^2 +3vw-1
\]
is a $2$-nct since it is irreducible.
For the irreducibility we refer the reader to Lemma~2.1 in \cite{GAGK}.
The area of $P_{\varphi_2}$ is $3/2$.

The polynomial
\[
\varphi_3 =  -\varphi_2(v-1) + v(w-1)^3
=  -1 + 6 v w - 4 v^2 w + v^3 w - 4 v w^2 + v^2 w^2 + v w^3
\]
is a $3$-nct since it is irreducible.
The area of $P_{\varphi_3}$ is $4$.

The polynomial
\begin{align*}
\varphi'_3 = & \ -\varphi_2(vw-1) - (v-1)^2(w-1)vw \\
= & \ -1 + 5 v w - 3 v^2 w + v^3 w - 2 v w^2 - v^2 w^2 + v^2 w^3
\end{align*}
is also a $3$-nct since it is irreducible.
The area of $P_{\varphi'_3}$ is $4$ if the characteristic of $k$ is not $2$.
(The Newton polygon $P_{\varphi'_3}$ depends on the characteristic of the base field. 
See Example~\ref{ExampleOfNct} (3).)

If the characteristic of $k$ is not $2$, the Newton polygons $P_{\varphi_i}$'s and  $P_{\varphi'_3}$ are as follows:

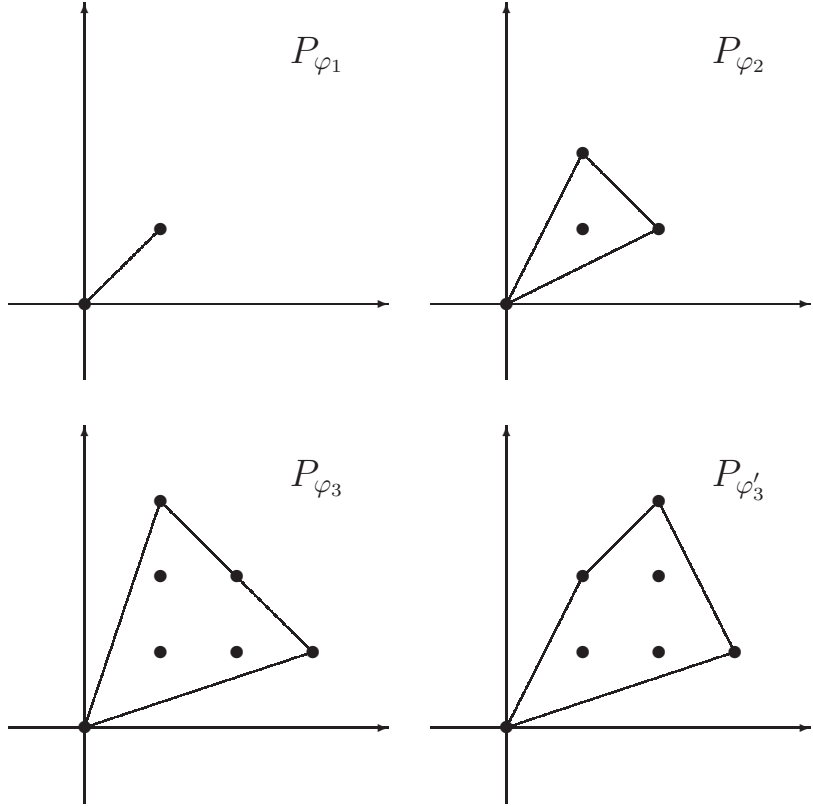
\[
{
\setlength\unitlength{1truecm}
\begin{picture}(5,5)(-1,-3)
\put(2.7,1.2){\mbox{\large $P_{\varphi_1}$}}
\qbezier (0,-2) (0.5,-1.5) (1,-1)
\put(-1,-2){\vector(1,0){5}}
\put(0,-3){\vector(0,1){5}}
\put(-0.1,-2.1){$\bullet$}
\put(0.9,-1.1){$\bullet$}
  \end{picture}
} \ \ \ \ 
{
\setlength\unitlength{1truecm}
  \begin{picture}(5,5)(-1,-3)
  \put(2.7,1.2){\mbox{\large $P_{\varphi_2}$}}
\put(-1,-2){\vector(1,0){5}}
\put(0,-3){\vector(0,1){5}}
\qbezier (0,-2) (0.5,-1) (1,0)
\qbezier (0,-2) (1,-1.5) (2,-1)
\qbezier (1,0) (1.5,-0.5) (2,-1)
\put(-0.1,-2.1){$\bullet$}
\put(0.9,-1.1){$\bullet$}
\put(0.9,-0.1){$\bullet$}
\put(1.9,-1.1){$\bullet$}
  \end{picture}
} 
\]

\[
{
\setlength\unitlength{1truecm}
  \begin{picture}(5,5)(-1,-3)
    \put(2.7,1.2){\mbox{\large $P_{\varphi_3}$}}
\put(-1,-2){\vector(1,0){5}}
\put(0,-3){\vector(0,1){5}}
\qbezier (0,-2) (0.5,-0.5) (1,1)
\qbezier (0,-2) (1.5,-1.5) (3,-1)
\qbezier (1,1) (2,0) (3,-1)
\put(-0.1,-2.1){$\bullet$}
\put(0.9,-1.1){$\bullet$}
\put(0.9,-0.1){$\bullet$}
\put(0.9,0.9){$\bullet$}
\put(1.9,-1.1){$\bullet$}
\put(1.9,-0.1){$\bullet$}
\put(2.9,-1.1){$\bullet$}
  \end{picture}
} \ \ \ \ 
{
\setlength\unitlength{1truecm}
  \begin{picture}(5,5)(-1,-3)
      \put(2.7,1.2){\mbox{\large $P_{\varphi'_3}$}}
\put(-1,-2){\vector(1,0){5}}
\put(0,-3){\vector(0,1){5}}
\qbezier (0,-2) (0.5,-1) (1,0)
\qbezier (0,-2) (1.5,-1.5) (3,-1)
\qbezier  (1,0) (1.5,0.5) (2,1)
\qbezier (2,1) (2.5,0)  (3,-1)
\put(-0.1,-2.1){$\bullet$}
\put(0.9,-1.1){$\bullet$}
\put(0.9,-0.1){$\bullet$}
\put(1.9,0.9){$\bullet$}
\put(1.9,-1.1){$\bullet$}
\put(1.9,-0.1){$\bullet$}
\put(2.9,-1.1){$\bullet$}
  \end{picture}
} 
\]

We can inductively construct an $r$-nct $\varphi_r$ by
\[
\varphi_r = -\varphi_{r-1}(v-1) + (-1)^{r-1}v(w-1)^r 
\]
for $r \ge 2$ as in Gonz\'alezAnaya-Gonz\'alez-Karu~\cite{GAGK}.
\end{rm}
\end{Example}

Negative curves are deeply related to $r$-ncts as follows:

\begin{Proposition}\label{Prop3.2}
Let $S_{a,b,c}$, $\bmp_{a,b,c}$ be the ring and the ideal as in Definition~\ref{Def1}.
Let $d_1$ and $r_1$ be positive integers.
Suppose that $f \in (S_{a,b,c})_{d_1}$ corresponds to $\varphi \lambda^{d_1} \in \left( \bigoplus_{(\alpha, \beta) \in d_1 P_{a,b,c} \cap \bZ^2} kv^\alpha w^\beta \right) \lambda^{d_1}$ under the identification (\ref{iden}).\footnote{
Remark that $N_\varphi$ is contained in $d_1 P_{a,b,c}$.}

Then $f$ is irreducible in $S_{a,b,c}$ and contained in $[{\bmp_{a,b,c}}^{(r_1)}]_{d_1}$ with 
$\frac{d_1}{r_1} < \sqrt{abc}$ (that is, $f \in [{\bmp_{a,b,c}}^{(r_1)}]_{d_1}$ is a negative curve) if and only if the following four conditions are satisfied:
\begin{enumerate}
\item[(1)]
$\varphi$ is irreducible in $A = k[v^{\pm 1}, w^{\pm 1}]$.
\item[(2)]
There exists an element of $N_\varphi$ on each edge of $d_1 P_{a,b,c}$.
\item[(3)]
$\varphi \in (v-1,w-1)^{r_1}A$.
\item[(4)]
$2|d_1 P_{a,b,c}| < {r_1}^2$.
\end{enumerate}
When this is the case, $\varphi$ is an $r_1$-nct.
\end{Proposition}

\Proof
Recall that $|P_{a,b,c}| = 1/2abc$ as in Remark~\ref{Area}.
Therefore the condition (4) is equivalent to $\frac{d_1}{r_1} < \sqrt{abc}$.
The condition (3) is equivalent to $f \in [{\bmp_{a,b,c}}^{(r_1)}]_{d_1}$
by (\ref{nth}).

Recall $r_1 > 0$.
If either the condition (3) or $f \in [{\bmp_{a,b,c}}^{(r_1)}]_{d_1}$ is satisfied,
$N_\varphi$ contains at least two elements.
Then, by Lemma~\ref{PrimeElementOfE}, 
$\varphi \lambda^{d_1}$ is irreducible in $S_{a,b,c} = E(P_{a,b,c}, \lambda)$ if and only if
both (1) and (2) are satisfied.
\qed

\vspace{2mm}

Here we state basic properties of $r$-ncts.

\begin{Proposition}\label{Fact3.3}
Let $k$ be a field and $A = k[v^{\pm 1}, w^{\pm 1}]$ be a Laurent polynomial ring with two variables $v$, $w$.
We put $\bfp=(v-1,w-1)A$.
Let $r$ be a positive integer.
\begin{enumerate}
\item
A unit element of $A$ is of the form $cv^\alpha w^\beta$, where $c \in k^\times$, $\alpha, \beta \in \bZ$.
If $u$ is a unit of $A$ and $\varphi$ is an $r$-nct,
then $u\varphi$ is also an $r$-nct.

Let $\xi: A \rightarrow A$ be a $k$-isomorphism such that $\xi(\bfp) = \bfp$.
Then there exists $(a_{ij}) \in {\rm GL}(2, \bZ)$ such that $\xi(v) = v^{a_{11}}w^{a_{12}}$, $\xi(w) = v^{a_{21}}w^{a_{22}}$.
If $\varphi$ is an $r$-nct, then so is $\xi(\varphi)$.
\item
Assume that $\varphi$ is an $r$-nct over $k$.
\begin{enumerate}
\item
We have
\begin{equation}\label{1dimensional}
\bfp^r \cap \left( \bigoplus_{(\alpha, \beta) \in P_\varphi \cap \bZ^2} kv^\alpha w^\beta \right) = k \varphi .
\end{equation}
\item
The number of $P_\varphi \cap \bZ^2$ is at most
$\frac{r(r+1)}{2} + 1$.
\item
If $\zeta$ is a reducible element contained in $\bfp^r$,
then $P_{\varphi}$ does not contain $N_\zeta$.
In particular, if $r \ge 2$, $P_{\varphi}$ does not contain $r+1$ points on a line.
\end{enumerate}
\item
Assume that $\varphi$ is an $r$-nct over $k$.
Then there exists an element $c \in k^\times$ such that all the coefficients of $c\varphi$ is in the prime field of $k$.
\item
If $\varphi$ is an $r$-nct over $k$, then $\varphi \not\in \bfp^{r+1}$.
\item
Assume that $\varphi$ is an $r$-nct over $k$.
Let $L/k$ be a field extension.
Then $\varphi$ is also an $r$-nct over $L$.
\item
Let $NCT_r$ be the set of $r$-ncts over $k$.
Consider the equivalence relation $\sim$ on $NCT_r$ generated by $\varphi \sim u\varphi$ and $\varphi \sim \xi(\varphi)$ as in (1). 
Then the quotient set $NCT_r/ \! \mathop{\sim}$ is a non-empty finite set for each $r$.
\end{enumerate}
\end{Proposition}

\Proof
It is easy to show (1).  
We omit a proof of it.

We shall prove (2) (a).
Assume the contrary.
Let $\varphi'$ be an element contained in the left-hand side but not in the right-hand one.
Since $\varphi \lambda$ is a prime element of $E(P_\varphi, \lambda)$ by Lemma~\ref{PrimeElementOfE},
$\varphi \lambda$, $\varphi' \lambda$ is an $E(P_\varphi, \lambda)$-regular sequence.
Take a homogeneous element $h \in E(P_\varphi, \lambda)$ such that
$\varphi \lambda$, $\varphi' \lambda$, $h$ is a homogeneous system of parameters
of $E(P_\varphi, \lambda)$.
Since $E(P_\varphi, \lambda)$ is a $3$-dimensional graded Cohen-Macaulay ring, we have
\begin{align*}
& \ \ell_{E(P_\varphi, \lambda)} \left( E(P_\varphi, \lambda) /(\varphi \lambda, \varphi' \lambda, h) \right)
= e\left( (h), E(P_\varphi, \lambda) /(\varphi \lambda, \varphi' \lambda) \right) \\
\ge & \ \ell_{E(P_\varphi, \lambda)_{\bmp_\varphi}} 
\left( E(P_\varphi, \lambda)_{\bmp_\varphi} /(\varphi \lambda, \varphi' \lambda) E(P_\varphi, \lambda)_{\bmp_\varphi} \right)
\cdot e\left( (h), E(P_\varphi, \lambda) /{\bmp_\varphi} \right) \\
\ge & \ r^2 (\deg h) ,
\end{align*}
where $\ell$ is the length, $e$ is the multiplicity and $\bmp_\varphi = E(P_\varphi, \lambda) \cap \bfp A[\lambda]$.
Here we remark that the first inequality comes from the additive formula of multiplicities,
and the second one comes from $\varphi \lambda, \varphi' \lambda \in \bmp_\varphi^r$
and $e((h), E(P_\varphi,\lambda)/\bmp_\varphi)) = e((\lambda^{\deg h}), k[\lambda]) = \deg h$.
On the other hand, consider the Poincare series
\[
P(E(P_\varphi, \lambda), s) = \sum_{n \ge 0} \dim_k E(P_\varphi, \lambda)_n s^n = \frac{f(s)}{(1-s)^3} ,
\]
where we remark that $E(P_\varphi, \lambda)$ is generated by elements of degree $1$ over $k$.
Here $f(1)$ is equal to the multiplicity of $E(P_\varphi, \lambda)$, which is $2|P_\varphi|$ by a theorem of Ehrhart (e.g. see Part~II of \cite{Hibi}).
Hence we have
\[
P(E(P_\varphi, \lambda)/(\varphi \lambda, \varphi' \lambda, h), s) =
\frac{f(s)(1-s)^2(1-s^{\deg h})}{(1-s)^3} 
\]
since $\varphi \lambda$, $\varphi' \lambda$, $h$ is a regular sequence.
Substituting $1$ for $s$, we have 
\[
\ell_{E(P_\varphi, \lambda)} \left( E(P_\varphi, \lambda) /(\varphi \lambda, \varphi' \lambda, h) \right)
= 2|P_\varphi| (\deg h) < r^2 (\deg h) .
\]
It is a contradiction.
Thus the equality (\ref{1dimensional}) is proved.\footnote{We can also prove (\ref{1dimensional}) using intersection theory on the blow-up of $\proj(E(P_\varphi, \lambda))$ at $(1,1)$ in the torus.}

The left-hand side of  (\ref{1dimensional})  is defined in $\bigoplus_{(\alpha, \beta) \in P_\varphi \cap \bZ^2} kv^\alpha w^\beta$ by $\frac{r(r+1)}{2}$ linear equations.
Therefore (b) follows from (a).

We shall prove (c).
If $\zeta$ is a reducible element contained in $\bfp^r$,
then $k \varphi$ does not contain $\zeta$ and $P_{\varphi}$ does not contain $N_\zeta$.
Assume that $P_{\varphi}$ contains $r+1$ points on a line.
Replacing $\varphi$ by some $u \xi(\varphi)$ as in (1),
we may assume that $P_{\varphi}$ contains $(0,0)$, $(1,0)$, \ldots, $(r,0)$.
(Here recall that $P_\varphi$ is a convex polygon.)
Then we have $(v-1)^r \in \bfp^r$ and $N_{(v-1)^r} \subset P_\varphi$.
It is a contradiction since $(v-1)^r$ is reducible for $r \ge 2$.

We shall prove (3).
Let $F$ be the prime field of $k$.
The assertion immediately comes from
\begin{align*}
& \ \left( (v-1,w-1)^r F[v^{\pm 1}, w^{\pm 1}] \cap \left( \bigoplus_{(\alpha, \beta) \in P_\varphi \cap \bZ^2} Fv^\alpha w^\beta \right) \right) \otimes_Fk \\
= & \ (v-1,w-1)^r k[v^{\pm 1}, w^{\pm 1}] \cap \left( \bigoplus_{(\alpha, \beta) \in P_\varphi \cap \bZ^2} kv^\alpha w^\beta \right) = k \varphi .
\end{align*}

Next we shall prove (4).
Suppose that an $r$-nct satisfies $\varphi \in \bfp^{r+1}$. 
Multiplying $\varphi$ by a unit of $A$, we may assume that 
$\varphi$ is a Laurent polynomial over the prime field and  
the origin is in $N_\varphi$.
Assume that the characteristic of $k$ is a prime number $p$.
(In the case of characteristic $0$, we can prove the assertion easier.)
If both $\alpha$ and $\beta$ were divided by $p$ for any $(\alpha,\beta) \in N_\varphi$,
then $\varphi$ would be reducible since $\varphi$ is a Laurent polynomial over the prime field of characteristic $p$.
Therefore we may assume that there exists $(\alpha', \beta') \in N_\varphi$ such that $p \not| \ \alpha'$.
Then we have
\[
0 \neq v \frac{\partial \varphi}{\partial v} \in \bfp^r \cap \left( \bigoplus_{(\alpha, \beta) \in P_\varphi \cap \bZ^2} kv^\alpha w^\beta \right)  = k \varphi .
\]
Here remark that the constant term of $v \frac{\partial \varphi}{\partial v}$ is $0$.
It is a contradiction since $\varphi$ and $v \frac{\partial \varphi}{\partial v}$ are linearly independent over $k$. 

We shall prove (5).
We have only to show that $\varphi$ is irreducible in $L[v^{\pm 1}, w^{\pm 1}]$.
Assume the contrary and suppose $\varphi = \psi_1\psi_2$ in $L[v^{\pm 1}, w^{\pm 1}]$.
Then we have 
\[
P_\varphi = P_{\psi_1} + P_{\psi_2} ,
\]
where the right-hand side is the Minkowski sum, that is $\{ \bma + \bmb \mid 
\bma \in P_{\psi_1}, \bmb \in P_{\psi_2} \}$.
Suppose
\begin{align*}
& \psi_1 \in 
(v-1,w-1)^{r_1} L[v^{\pm 1}, w^{\pm 1}] 
\setminus (v-1,w-1)^{r_1+1} L[v^{\pm 1}, w^{\pm 1}] , \\
& \psi_2 \in (v-1,w-1)^{r_2} L[v^{\pm 1}, w^{\pm 1}] \setminus (v-1,w-1)^{r_2+1} L[v^{\pm 1}, w^{\pm 1}] .
\end{align*}
Then, by (4), we have $r = r_1+r_2$.
Since $r^2 > 2|P_\varphi|$,  we have
\[
\frac{r_1}{\sqrt{2}} + \frac{r_2}{\sqrt{2}} = \frac{r}{\sqrt{2}}  > \sqrt{|P_\varphi|} \ge \sqrt{|P_{\psi_1}|} + \sqrt{|P_{\psi_2}|} .
\]
Here the second inequality is called the Brunn-Minkowski inequality.
Therefore either $r_1^2 > 2|P_{\psi_1}|$ or $r_2^2 > 2|P_{\psi_2}|$ is satisfied.
Hence we know that some irreducible divisor $\varphi'$ of $\varphi$ is an $r'$-nct for some $r'$.
By (3), we may assume that $\varphi'$ is a Laurent polynomial over the prime field.
It is a contradiction since $\varphi$ is irreducible over $k$.

We shall prove (6).
Since $NCT_r$ contains $\varphi_r$ in Example~\ref{123-nct},
the set $NCT_r/\sim$ is not empty.
Any $1$-nct is equivalent to $v-1$.
Suppose $r \ge 2$.
Let $\varphi$ be an $r$-nct.
If $|P_{\varphi}| =0$,
then all the points of $N_\varphi$ are on a line and we may assume that $\varphi$ is a polynomial in $v$.
Then $(v-1)^r$ divides $\varphi$.
It contradicts to the irreducibility of $\varphi$.
Therefore we have $|P_{\varphi}| >0$.
Let $\Omega$ be the convex hull of the four points $(0, 0)$, $(\sqrt{2}r^2, 0)$, $((\sqrt{2}+1)r^2, r^2)$, $(0, r^2)$.
Let $P$ be an integral convex polygon (a convex hull of a finite subset of $\bZ^2$ with positive area).
Then, by (2), $r$-ncts $\varphi$ with $P_\varphi = P$ are equivalent to each others.
Since $\Omega$ is bounded, there exist only finitely many integral convex polygons contained in $\Omega$.
Therefore there exists a finite subset $F$ of $NCT_r$ such that 
any $r$-nct $\varphi$ with $P_\varphi \subset \Omega$ is equivalent to 
one of $F$.
It is enough to prove that any integral convex polygon with area less than $r^2/2$ is contained in $\Omega$ after an affine transformation $\xi$ satisfying $\xi(\bZ^2) = \bZ^2$.
Let $P$ be an  integral convex polygon with area less than $r^2/2$.
By an affine transformation preserving the lattice, we may assume that 
$(0,0)$ and $(\alpha_1, 0)$ are adjacent vertices of $P$, where $\alpha_1$ is a positive integer.
Since $|P| < r^2/2$, we may assume that any point $(\alpha,\beta)$ in $P$ satisfies $0 \le \beta < r^2$.
Let $(\alpha_2, \beta_2)$ be the other vertex adjacent to $(0,0)$.
By a linear transformation of the form $\left( \begin{array}{cc} 1 & a \\ 0 & 1\end{array} \right)$ for some $a \in \bZ$, 
we may assume $\beta_2> \alpha_2\ge 0$.
Then any point $(\alpha,\beta)$ in $P$ satisfies $0 \le \alpha$.
Suppose that $P$ is not contained in $\Omega$.
Take $(\alpha_3,\beta_3) \in P \setminus \Omega$.
Let $\ell$ be the line through $(0,0)$ and $(\alpha_2, \beta_2)$.
The distance of the line $\ell$ and the point $(\alpha_3,\beta_3)$ is 
bigger than that of $\ell$ and $(\sqrt{2}r^2,0)$.
Therefore
\[
|P| > \ \mbox{(the area of the triangle $(0,0)$, $(\alpha_2,\beta_2)$, $(\sqrt{2}r^2,0)$)} 
=   \frac{\beta_2 \sqrt{2}r^2}{2} > \frac{r^2}{2} .
\]
It is a contradiction. Therefore we have $P \subset \Omega$.
\qed

\begin{Example}\label{ExampleOfNct}
\begin{rm}
\begin{enumerate}
\item
The quotient set $NCR_1/ \! \mathop{\sim}$ consists of the equivalence class of $\varphi_1$
in Example~\ref{123-nct}.
\item
The quotient set $NCR_2/ \! \mathop{\sim}$ consists of the equivalence class of $\varphi_2$ in Example~\ref{123-nct}.
In fact, suppose that $\varphi$ is a $2$-nct.
By Proposition~\ref{Fact3.3} (2) (c), $P_\varphi$ does not contain three points on a line. 
As in the proof of Proposition~\ref{Fact3.3} (6), we may assume that $P_\varphi$ has three successive vertices $(\alpha_2, \beta_2)$, $(0,0)$, $(1,0)$ where $0 \le \alpha_2 < \beta_2 < 4$.
If $(\alpha_2, \beta_2)=(0,1)$, then $P_\varphi$ contains $(1,1)$.
Putting $\zeta = (v-1)(w-1) \in \bfp^2$, $P_\varphi$ contains $N_\zeta$.
It contradicts to Proposition~\ref{Fact3.3} (2) (c).
Then we know $(\alpha_2, \beta_2) = (2,3)$.
In this case, $\varphi$ is equivalent to $\varphi_2$.

In the case of $r = 1, 2$, $NCR_r/ \! \mathop{\sim}$ consists of the equivalence class of $\varphi_r$ and
$P_{\varphi_r}$ is independent of $k$.
Therefore,  for $r= 1, 2$, the existence of negative curves contained in $[{\bmp_{a,b,c}}^{(r)}]_d$ does not depend on the base field $k$ (see Proposition~\ref{Prop3.2}).
\item
One can prove that
the quotient set $NCT_3/ \! \mathop{\sim}$ consists of two elements, which are represented by 
$\varphi_3$ and $\varphi'_3$ in Example~\ref{123-nct}, respectively.
The polygon $P_{\varphi_3}$ is the triangle with vertices $(0,0)$, $(3,1)$, $(1,3)$.

If the characteristic of $k$ is not $2$, then $P_{\varphi'_3}$ is the tetragon as in Example~\ref{123-nct}.

Assume that the characteristic of $k$ is $2$.
Since the coefficient of $vw^2$ in $\varphi'_3$ is $-2$, 
$\varphi'_3$ is a $3$-nct such that $P_{\varphi'_3}$ is the triangle with vertices
$(0,0)$, $(3,1)$ and $(2,3)$.
In the case where $(a,b,c)=(9,10,13)$ and ${\rm ch}(k) = 2$, there exists a negative curve $f \in [{\bmp_{9,10,13}}^{(3)}]_{100}$ 
corresponding to a $3$-nct which is equivalent to $\varphi'_3$  as in the proof of Theorem~4.1 in \cite{MG}.
There does not exist a negative curve in ${\bmp_{9,10,13}}^{(3)}$ if ${\rm ch}(k) \neq 2$.
\end{enumerate}
\end{rm}
\end{Example}

\section{The Cox ring of the blow-up of $X_{\Delta_{\varphi}}$ and 
the symbolic Rees ring of the Ehrhart ring $E(P_{\varphi}, \lambda)$}
\label{sect4}

Let $f \in [{\bmp_{a,b,c}}^{(r_1)}]_{d_1}$ be a negative curve and $\varphi$ be an $r_1$-nct
such that $f = \varphi \lambda^{d_1}$ as in Proposition~\ref{Prop3.2}.
Then we have 
\[
|P_\varphi| \le |d_1 P_{a,b,c}| < r^2/2 .
\]
From the results of Gonz\'alezAnaya-Gonz\'alez-Karu~\cite{GAGK},
it can be inferred that $R'_s(\bmp_{a,b,c})$ tends to be finitely generated when $|d_1 P_{a,b,c}|$ is close to $|P_\varphi|$, and
$R'_s(\bmp_{a,b,c})$ tends to be infinitely generated when $|d_1 P_{a,b,c}|$ is close to $r^2/2$.
Therefore it is natural to ask:
\[
\mbox{Is $R'_s(\bmp_\varphi)$ a Noetherian ring?}
\]
Here $\bmp_\varphi$ is a prime ideal of the Ehrhart ring $E(P_\varphi,\lambda) \subset k[v^{\pm 1}, w^{\pm 1}][\lambda]$ defined by $\bmp_\varphi = E(P_\varphi,\lambda) \cap (v-1,w-1)k[v^{\pm 1}, w^{\pm 1}][\lambda]$.

We put $X_{\Delta_{a,b,c}} = \proj(E(P_{a,b,c}, \lambda))$ and 
$X_{\Delta_\varphi} = \proj(E(P_\varphi, \lambda))$, where $\Delta_{a,b,c}$ and $\Delta_\varphi$ are the fans corresponding to the toric varieties $X_{\Delta_{a,b,c}}$ and $X_{\Delta_\varphi}$, respectively.
Let $Y_{\Delta_{a,b,c}}$ (resp.\ $Y_{\Delta_\varphi}$) be the blow-up of $X_{\Delta_{a,b,c}}$ (resp.\ $X_{\Delta_\varphi}$) at the point $(1,1)$ in the torus $\spec k[v^{\pm 1}, w^{\pm 1}]$.
Another reason why we are studying $R'_s(\bmp_\varphi)$ is that
$Y_{\Delta_\varphi}$ contains a negative curve that is birational to the negative curve in $Y_{\Delta_{a,b,c}}$.
We are interested in the birational class of the negative curve in $Y_{\Delta_{a,b,c}}$.
(The author does not know any example that  the negative curve in $Y_{\Delta_{a,b,c}}$ is not rational.)

\vspace{2mm}

We shall prove Theorem~\ref{Thm3.6} in this section.
Before proving this theorem, we shall show the following lemma:

\begin{Lemma}\label{Lemma8.9}
With notation as in Theorem~\ref{Thm3.6}, the following conditions are equivalent:
\begin{enumerate}
\item[(i)]
$H^0(Y_{\Delta_\varphi}, {\cal O}_{Y_{\Delta_\varphi}}(K_{Y_{\Delta_\varphi}} + n C_\varphi)) = 0$ for any $n > 0$,
\item[(ii)]
$H^0(Y_{\Delta_\varphi}, {\cal O}_{Y_{\Delta_\varphi}}(K_{Y_{\Delta_\varphi}} +  C_\varphi)) = 0$,
\item[(iii)]
$C_\varphi . (K_{Y_{\Delta_\varphi}} + C_\varphi) = -2$,
\item[(iv)]
$C_\varphi . (K_{Y_{\Delta_\varphi}} + C_\varphi) < 0$,
\item[(v)]
$I_\varphi = \frac{r(r-1)}{2}$,
\item[(vi)]
$I_\varphi \le \frac{r(r-1)}{2}$.
\end{enumerate}

If $C_\varphi \simeq \bP_k^1$, then {\rm (ii)} is satisfied.
In the case where $k$ is algebraically closed, the converse is also true.
\end{Lemma}

\Proof
${\rm (i)} \Rightarrow {\rm (ii)}$, ${\rm (iii)} \Rightarrow {\rm (iv)}$, ${\rm (v)} \Rightarrow {\rm (vi)}$ are obvious.

We shall prove ${\rm (iv)} \Rightarrow {\rm (ii)}$.
Assume that there exists an effective Weil divisor $D$ on $Y_{\Delta_\varphi}$
that is linearly equivalent to $K_{Y_{\Delta_\varphi}}+C_\varphi$.
By (iv), there exists an effective Weil divisor $D'$ such that $D = C_\varphi+ D'$.
Let $\pi_\varphi: Y_{\Delta_\varphi} \rightarrow X_{\Delta_\varphi}$ be the blow-up at the point $(1,1)$.
Since
\[
K_{X_{\Delta_\varphi}} + (\pi_\varphi)_*(C_\varphi) \sim (\pi_\varphi)_*(D) = (\pi_\varphi)_*(C_\varphi) + (\pi_\varphi)_*(D') ,
\]
we have
\[
(\pi_\varphi)_*(D') - K_{X_{\Delta_\varphi}} \sim 0 .
\]
It contradicts to the fact that the left-hand side is a non-zero effective divisor.

We shall prove ${\rm (ii)} \Rightarrow {\rm (vi)}$.
We have 
\[
0 = H^0(Y_{\Delta_\varphi}, {\cal O}_{Y_{\Delta_\varphi}}(K_{Y_{\Delta_\varphi}} +  C_\varphi))
= \bfp^{r-1} \cap \left( \bigoplus_{(\alpha,\beta) \in {P_\varphi}^\circ \cap \bZ^2}
kv^\alpha w^\beta \right) ,
\]
where ${P_\varphi}^\circ$ is the interior of $P_\varphi$ and $\bfp = (v-1,w-1)k[v^{\pm 1}, w^{\pm 1}]$.
Since $\bfp^{r-1}$ is defined by $\frac{r(r-1)}{2}$ linear equations in $k[v^{\pm 1}, w^{\pm 1}]$,
the number of ${P_\varphi}^\circ \cap \bZ^2$ is less than or equal to $\frac{r(r-1)}{2}$.

Next we shall prove ${\rm (iii)} \Leftrightarrow {\rm (v)}$ and ${\rm (vi)} \Rightarrow {\rm (iv)}$.
By Pick's theorem, we have
\[
2|P_\varphi| = B_\varphi + 2I_\varphi -2 .
\]
Since $C_\varphi^2 = 2|P_\varphi| - r^2$ and $C_\varphi . (-K_{Y_{\Delta_\varphi}}) = B_\varphi - r$, we have
\[
C_\varphi . (K_{Y_{\Delta_\varphi}} + C_\varphi)  + 2 
=  - B_\varphi + r +2|P_\varphi| - r^2 +2
=  2I_\varphi - r(r-1) .
\]
It is easy to check ${\rm (iii)} \Leftrightarrow {\rm (v)}$ and ${\rm (vi)} \Rightarrow {\rm (iv)}$.

We shall prove ${\rm (ii)} \Rightarrow {\rm (i)}$.
As we have already seen, (ii) is equivalent to (iv).
Therefore we have $C_\varphi . (K_{X_{\Delta_\varphi}} + nC_\varphi) < 0$ for any $n > 0$.
Since any effective divisor that is linearly equivalent to $K_{Y_{\Delta_\varphi}} + nC_\varphi$ has $C_\varphi$ as a component, the multiplication by $C_\varphi$ 
\[
H^0(Y_{\Delta_\varphi}, {\cal O}_{Y_{\Delta_\varphi}}(K_{Y_{\Delta_\varphi}} + (n-1) C_\varphi)) \rightarrow H^0(Y_{\Delta_\varphi}, {\cal O}_{Y_{\Delta_\varphi}}(K_{Y_{\Delta_\varphi}} + n C_\varphi))
\]
is bijective.
Therefore (ii) implies (i).

Next we shall prove ${\rm (vi)} \Rightarrow {\rm (v)}$.
We may assume that $k$ is algebraically closed.
Remark that 
\begin{equation}\label{P1}
h^0({\cal O}_{Y_{\Delta_\varphi}}( K_{Y_{\Delta_\varphi}} + C_\varphi))
= h^2({\cal O}_{Y_{\Delta_\varphi}}( - C_\varphi)) = h^1({\cal O}_{C_\varphi}).
\end{equation}
Since (vi) is equivalent to (ii), we have $h^1({\cal O}_{C_\varphi}) = 0$ and $C_\varphi \simeq \bP_k^1$.
By definition, $C_\varphi$ does not meet the singular points of $Y_{\Delta_\varphi}$. 
Then, by the adjunction formula, we have $\omega_{Y_{\Delta_\varphi}}( C_\varphi)|_{C_\varphi} = \omega_{C_\varphi}$. 
Since (iii) is satisfied, (v) holds.

If $C_\varphi \simeq \bP^1_k$, (ii) is satisfied by (\ref{P1}).
(ii) implies $C_\varphi \simeq \bP_k^1$ when $k$ is algebraically closed.
\qed

\vspace{2mm}

Now we start to prove Theorem~\ref{Thm3.6}.
$(1)\Rightarrow (2)$ is trivially true.

We shall prove
$(1)\Rightarrow (3)\Rightarrow  (4)\Rightarrow (5)\Rightarrow  (6)\Rightarrow  (10)$
for any integral convex polygon $P$ and the corresponding complete $2$-dimensional fan $\Delta$.

$(1)\Rightarrow (3)$ is a basic fact.
(Any nef and big $\bQ$-Cartier divisor $D$ on a normal projective surface satisfies $D^2 > 0$.)

We shall show $(3)\Rightarrow (4)$.
Recall that $-K_{X_\Delta}= D_1 + \cdots + D_s$ is a $\bQ$-Cartier divisor, where each $D_i$ is a toric prime divisor corresponding to each edge of $P$.
Take a positive integer $q$ such that $qP_{-K_{X_\Delta}}$ is an integral convex polygon
and $-qK_{X_\Delta}$ is a Cartier divisor on $X_\Delta$. 
By the Riemann-Roch theorem, we know that $\chi({\cal O}_{X_{\Delta}}(-nqK_{X_{\Delta}}))$ is a polynomial in $n$ of degree $2$ such that the coefficient of $n^2$ is $(-qK_{X_{\Delta}})^2/2$.
On the other hand, 
$h^0({\cal O}_{X_{\Delta}}(-nqK_{X_{\Delta}})) = {}^\#(nqP_{-K_{X_{\Delta}}} \cap \bZ^2)$ is a polynomial of degree $2$ for $n\ge 0$ (it is called the Ehrhart polynomial) and 
the coefficient of $n^2$ is $|qP_{-K_{X_{\Delta}}}|$.
Since
\[
\chi({\cal O}_{X_{\Delta}}(-nqK_{X_{\Delta}})) = 
h^0({\cal O}_{X_{\Delta}}(-nqK_{X_{\Delta}}))-h^1({\cal O}_{X_{\Delta}}(-nqK_{X_{\Delta}}))+h^2({\cal O}_{X_{\Delta}}(-nqK_{X_{\Delta}}))
\]
and $h^2({\cal O}_{X_{\Delta}}(-nqK_{X_{\Delta}})) = h^0({\cal O}_{X_{\Delta}}(K_{X_\Delta}+nqK_{X_{\Delta}})) = 0$ for $n \ge 0$,
we obtain 
\[
(-K_{X_{\Delta}})^2/2 \le |P_{-K_{X_{\Delta}}}| .
\]
Let $E$ be the exceptional divisor of the blow-up $Y_\Delta \rightarrow X_\Delta$
at $(1,1)$ in the torus.
Since $(-K_{Y_{\Delta}})^2 = (-K_{X_{\Delta}})^2 + E^2=(-K_{X_{\Delta}})^2-1$,
we obtain $|P_{-K_{X_{\Delta}}}| > \frac{1}{2}$.

We shall prove $(4)\Rightarrow (5)$.
We have
$h^0({\cal O}_{X_\Delta}( -nq K_{X_{\Delta}})) = {}^\#(nqP_{-K_{X_{\Delta}}} \cap \bZ^2) 
$.
By definition, $H^0(Y_\Delta, {\cal O}_{Y_\Delta}(-nq K_{Y_{\Delta}}))$ is a $k$-vector subspace of 
$H^0(X_\Delta,  {\cal O}_{X_\Delta}(-nq K_{X_{\Delta}}))$ defined by $nq(nq+1)/2$ linear equations.
Therefore we have
\[
{}^\#(nqP_{-K_{X_{\Delta}}} \cap \bZ^2) - \frac{nq(nq+1)}{2} \le h^0({\cal O}_{Y_\Delta}( -nq K_{Y_{\Delta}})) .
\]
Since ${}^\#(nqP_{-K_{X_{\Delta}}} \cap \bZ^2)$ is a polynomial of degree $2$ for $n\ge 0$ whose coefficient of $n^2$ is $|qP_{-K_{X_{\Delta}}}|$,
we know that $-K_{Y_{\Delta}}$ is big if $|P_{-K_{X_{\Delta}}}| > \frac{1}{2}$.

Next we shall prove $(5)\Rightarrow (6)$.
We may assume that $k$ is an algebraically closed field.

First we shall construct a refinement $\Delta'$ of $\Delta$ such that
\begin{itemize}
\item
$X_{\Delta'}$ is a smooth toric variety, and
\item
$P_{-K_{X_\Delta}}=P_{-K_{X_{\Delta'}}}$ .
\end{itemize}
Let 
\[
\{ \bRo \bma_1, \bRo \bma_2, \ldots, \bRo \bma_n \}
\]
be the set of the $1$-dimensional cones in $\Delta$.
We assume that each $\bma_i$ is the shortest integer vector in the cone $\bRo \bma_i$.
Assume that $\bma_1$, $\bma_2$, \ldots, $\bma_n$ are arranged counterclockwise around the origin.
We think that each $\bma_i$ is a row vector of length 2.
We shall construct $\Delta'$ using induction on 
\begin{equation}\label{bma}
\left|
{\rm det}\left( \begin{array}{l} \bma_1 \\ \bma_2 \end{array} \right) \times 
{\rm det}\left( \begin{array}{l} \bma_2 \\ \bma_3 \end{array} \right)  \times \cdots \times 
{\rm det}\left( \begin{array}{l} \bma_{n-1} \\ \bma_n \end{array} \right)  \times 
{\rm det}\left( \begin{array}{l} \bma_n \\ \bma_1 \end{array} \right)
\right| .
\end{equation}
Suppose that $X_\Delta$ is not smooth.
After a linear transformation in ${\rm SL}(2,\bZ)$,
we may assume that $\bma_1 = (1,0)$ and $\bma_2 = (a,b)$, where $b > a > 0$.
Here we put $\bmb = (1,1)$.
Then
\[
\left|
{\rm det}\left( \begin{array}{l} \bma_1 \\ \bmb \end{array} \right)  \times 
{\rm det}\left( \begin{array}{l} \bmb \\ \bma_2 \end{array} \right)  \times 
{\rm det}\left( \begin{array}{l} \bma_2 \\ \bma_3 \end{array} \right)  \times \cdots \times 
{\rm det}\left( \begin{array}{l} \bma_{n-1} \\ \bma_n \end{array} \right)  \times 
{\rm det}\left( \begin{array}{l} \bma_n \\ \bma_1 \end{array} \right)
\right| 
\]
is strictly less than (\ref{bma}).
Let $\bar{\Delta}$ be the complete fan in $\bR^2$ with $1$-dimensional cones
\[
\{ \bRo \bma_1, \bRo \bmb, \bRo \bma_2, \ldots, \bRo \bma_n \} .
\]
Here one can check $P_{-K_{X_\Delta}}=P_{-K_{X_{\bar{\Delta}}}}$.
Repeating this process, we can construct $\Delta'$ satisfying the required conditions.

Let $Y_{\Delta'}$ be the blow-up of $X_{\Delta'}$ at the point $(1,1)$ in the torus.
Then $Y_{\Delta'} \rightarrow Y_\Delta$ is a resolution of singularities.
Then, since $P_{-K_{X_\Delta}}=P_{-K_{X_{\Delta'}}}$, we have
\[
H^0(Y_\Delta, {\cal O}_{Y_\Delta}(-nK_{Y_\Delta})) = \bfp^n \cap 
\left( \bigoplus_{(\alpha,\beta) \in nP_{-K_{X_\Delta}} \cap \bZ^2} kv^\alpha w^\beta \right) 
= H^0(Y_{\Delta'}, {\cal O}_{Y_{\Delta'}}(-nK_{Y_{\Delta'}})), 
\]
where $\bfp = (v-1,w-1)k[v^{\pm 1},w^{\pm 1}]$.
Hence, if $-K_{Y_\Delta}$ is big, so is $-K_{Y_{\Delta'}}$.
Then, by Theorem~1 in Testa-VarillyAlvarado-Velasco~\cite{TVAV},
${\rm Cox}(Y_{\Delta'})$ is finitely generated.
Since $Y_{\Delta'} \rightarrow Y_\Delta$ is surjective,
${\rm Cox}(Y_{\Delta})$ is also finitely generated by Theorem~1.1 in Okawa~\cite{Okawa}.
We have completed the proof of $(5)\Rightarrow (6)$.

We shall prove $(6)\Rightarrow (10)$.
Let $E$ be the exceptional divisor of 
$Y_{\Delta}\rightarrow X_{\Delta}=\proj(E(P, \lambda))$.
Let $H$ be the pullback of ${\mathcal O}_{X_{\Delta}}(1)$.
We put $\bmp = E(P,\lambda) \cap \bfp k[v^{\pm 1}, w^{\pm 1}, \lambda]$.
Then the extended symbolic Rees ring $R'_s(\bmp)$ coincides with
the multisection ring 
\[
\bigoplus_{r, d \in \bZ} H^0(Y_{\Delta}, {\cal O}_{Y_{\Delta}}(dH -rE)) ,
\]
which is a pure subring of ${\rm Cox}(Y_{\Delta})$.
Therefore, if ${\rm Cox}(Y_{\Delta})$ is finitely generated,
so is $R'_s(\bmp)$.

In the rest of this proof, we consider an $r$-nct $\varphi$,
the integral convex polygon $P_\varphi$ and the corresponding fan $\Delta_\varphi$.

We shall prove $(5)\Rightarrow (8)$.
Suppose that there exists an effective divisor $D$ that is linearly equivalent to $K_{Y_{\Delta_\varphi}} + nC_\varphi$.
Then we have 
$
n C_\varphi \sim D -K_{Y_{\Delta_\varphi}} 
$ and $C_\varphi$ is big.
It contradicts to 
\[
{C_\varphi}^2 = 2|P_\varphi| - r^2 < 0 .
\]

$(2)\Rightarrow (7)$ follows from
$(-K_{Y_{\Delta_\varphi}}).C_\varphi = B_\varphi - r$.

We shall prove $(7) \Rightarrow (8)$.
Since $0 \le B_\varphi - r = (-K_{\Delta_\varphi}). C_\varphi$ and $C_\varphi^2< 0$,
(iv) in Lemma~\ref{Lemma8.9} is satisfied.
Therefore (i) in Lemma~\ref{Lemma8.9} is satisfied.

$(8)\Leftrightarrow (9)$ is nothing but ${\rm (i)} \Leftrightarrow {\rm (v)}$ in Lemma~\ref{Lemma8.9}.

$(11) \Rightarrow (9)$ and (b) follow from Lemma~\ref{Lemma8.9}.

In the rest of the proof, we prove $(8)\Rightarrow (10)$ and (c).

When we prove (c), 
we may assume that $k$ is a finite field.
Then we can prove (10) in the same way as Theorem~1 in Cutkosky~\cite{C}.

Next we shall prove $(8)\Rightarrow (10)$.
By (c), we may assume that $k$ is a field of characteristic $0$.
Furthermore, we may assume that $k$ is the field of complex numbers $\bC$.
Since 
\[
0=h^0({\cal O}_{Y_{\Delta_\varphi}}(K_{Y_{\Delta_\varphi}} + nC_\varphi)) = 
h^2({\cal O}_{Y_{\Delta_\varphi}}(-nC_\varphi)) = h^1({\cal O}_{nC_\varphi})
\]
for any $n > 0$, we can contract $C_\varphi$, that is, there exists a birational morphism $\xi:Y_{\Delta_\varphi} \rightarrow Z$ such that $\xi(C_\varphi)$ is a point, where $Z$ is a normal projective surface with at most rational singularity.
Here put $H= (\pi_\varphi)^*({\cal O}_{X_{\Delta_\varphi}}(1))$, where $\pi_\varphi: Y_{\Delta_\varphi} \rightarrow X_{\Delta_\varphi}$ is the blow-up at $(1,1)$ in the torus.
Let $E$ be the exceptional divisor of $\pi_\varphi$.
Let $i$ and $j$ be positive integers such that $(iH-jE).C_\varphi = 0$.
Then there exists $n > 0$ such that $n\xi_*(iH-jE)$ is a very ample Cartier divisor on $Z$.
Then one can prove that $n(iH-jE)$ is a semi-ample Cartier divisor on $Y_{\Delta_\varphi}$.
It is easy to verify that $R'_s(\bmp_\varphi)$ is Noetherian.
\qed

\begin{Example}\label{ExampleOfGGK}
\begin{rm}
Suppose that the characteristic of $k$ is $0$.
In Gonz\'alezAnaya-Gonz\'alez-Karu~\cite{GAGK}, they constructed two distinct $r$-ncts $\varphi_r$ and $\varphi'_r$ over $k$ for each $r \ge 3$.
Here $P_{\varphi_r}$ is a triangle with vertices $(-1,-1)$, $(r-1,0)$ and $(0,r-1)$.
The polygon $P_{\varphi'_r}$ is a tetragon with vertices $(-1,-1)$, $(r-1,0)$, $(\frac{r-1}{2},r-1)$, $(\frac{r-3}{2},r-2)$ in the case where $r$ is odd, and
$(-1,-1)$, $(r-1,0)$, $(\frac{r}{2},r-2)$, $(\frac{r-2}{2},r-1)$ in the case where $r$ is even.
Both $P_{\varphi_r}$ and $P_{\varphi'_r}$ contain $\frac{r(r+1)}{2} + 1$ lattice points.
Both $Y_{\Delta_{\varphi_r}}$ and $Y_{\Delta_{\varphi'_r}}$ satisfy the condition (1) in Theorem~\ref{Thm3.6}.

We shall give an outline of the proof of the above assertions for $\varphi'_r$.
(One can prove the same assertions for $\varphi_r$ in the same way.)
Consider the tetragon $P$ with four vertices as above.
Obviously it contains $\frac{r(r+1)}{2} + 1$ lattice points.
Put $\bfp = (v-1,w-1)k[v^{\pm 1}, w^{\pm 1}]$.
Since
\begin{equation}\label{pr}
\bfp^r \cap \left( \bigoplus_{(\alpha,\beta) \in P \cap \bZ^2} kv^\alpha w^\beta \right) 
\end{equation}
is defined by $\frac{r(r+1)}{2}$ linear equations in $\left( \bigoplus_{(\alpha,\beta) \in P \cap \bZ^2} kv^\alpha w^\beta \right)$, 
(\ref{pr}) is not $0$.
Let $\varphi'_r$ be a non-zero element in (\ref{pr}).
Using Lemma~\ref{EU} below, we have 
\[
\bfp^r \cap \left( \bigoplus_{
\substack{ 
(\alpha,\beta) \in P \cap \bZ^2
  \\
(\alpha,\beta) \neq (-1,-1)
}}  
kv^\alpha w^\beta \right) = 0 .
\]
Therefore the coefficient of $v^{-1}w^{-1}$ in $\varphi'_r$ is not zero.
In the same way, we know that the coefficients of monomials corresponding to 
the vertices of $P$ are not zero.
Thus we obtain $P_{\varphi'_r} = P$.
By Lemma~2.1 in \cite{GAGK}, we know $\varphi'_r$ is irreducible.
Since $|P| = \frac{r^2-1}{2}$,
we know that $\varphi'_r$ is an $r$-nct.

Next we shall prove that $-K_{Y_{\Delta_{\varphi'_r}}}$ is nef and big.
It is easy to see that $|P_{-K_{X_{\Delta_{\varphi'_r}}}}| > \frac{1}{2}$.
Therefore $-K_{Y_{\Delta_{\varphi'_r}}}$ is big by Theorem~\ref{Thm3.6}.
Let $V$ be the closure of 
\[
\spec(k[v^{\pm 1}, w^{\pm 1}]/(w-1)) 
\]
in $X_{\Delta_{\varphi'_r}}$.
Let $D$ be the toric prime divisor on $X_{\Delta_{\varphi'_r}}$ corresponding to 
the bottom edge of $P_{\varphi'_r}$.
Let $\tilde{V}$ and $\tilde{D}$ be the proper transforms of $V$ and $D$, respectively.
Then we know $-K_{Y_{\Delta_{\varphi'_r}}}$ is linearly equivalent to $\tilde{V}+\tilde{D}$.
Since $V^2 > 1$ and $D^2 > 0$,
we obtain $\tilde{V}^2>0$ and $\tilde{D}^2>0$.
Thus we know $-K_{Y_{\Delta_{\varphi'_r}}}$ is nef.
\end{rm}
\end{Example}

\begin{Lemma}\label{EU}
Let $k$ be a field of characteristic $0$ and $n$ be a positive integer.
For a subset $U$ of $\bZ^2$, we put
\[
k\seisei{U} = \{ \sum_{(\alpha,\beta)\in U} c_{(\alpha,\beta)}v^\alpha w^\beta \mid c_{(\alpha,\beta)} \in k \} .
\]
Let $L$ be a line in $\bR^2$ such that ${}^\#(L \cap U) = n$.
Put $U' = U \setminus (L \cap U)$.
Then we have an isomorphism of $k$-vector spaces
\[
k\seisei{U} \cap (v-1,w-1)^nk[v^{\pm 1}, w^{\pm 1}] \simeq k\seisei{U'} \cap (v-1,w-1)^{n-1}k[v^{\pm 1}, w^{\pm 1}] .
\]
\end{Lemma}

We omit a proof of Lemma~\ref{EU}.
We can prove this lemma in the same way as the proof of Lemma~4.5 in \cite{KN}.

\begin{Remark}\label{Ex3.7}
\begin{rm}
Let $\varphi$ be an $r$-nct.
By Proposition~\ref{Fact3.3} (2) (b), we have ${}^\#(P_\varphi \cap \bZ^2) \le \frac{r(r+1)}{2} + 1$.
In many cases ${}^\#(P_\varphi \cap \bZ^2) = \frac{r(r+1)}{2} + 1$ is satisfied.
Here we show that it is equivalent to $B_\varphi = r+1$.
If $B_\varphi = r+1$ is satisfied, we know $I_\varphi = \frac{r(r-1)}{2}$ by Theorem~\ref{Thm3.6}.
Thus we obtain 
\[
{}^\#(P_\varphi \cap \bZ^2) = B_\varphi + I_\varphi = \frac{r(r+1)}{2} + 1 .
\]
Conversely assume $B_\varphi + I_\varphi = \frac{r(r+1)}{2} + 1$.
Since
\[
0 < r^2-2|P_\varphi| = r^2 - 2I_\varphi - B_\varphi + 2 = 
r^2 - 2(I_\varphi + B_\varphi) + B_\varphi + 2 = B_\varphi -r 
\]
by Pick's theorem, we know $I_\varphi = \frac{r(r-1)}{2}$ by Theorem~\ref{Thm3.6}.
Therefore we have $B_\varphi = {}^\#(P_\varphi \cap \bZ^2) - I_\varphi = r+1$.

In particular, the condition (7) in Theorem~\ref{Thm3.6} is satisfied if ${}^\#(P_\varphi \cap \bZ^2) = \frac{r(r+1)}{2} + 1$.
\end{rm}
\end{Remark}

\begin{Remark}\label{SMC}
\begin{rm}
Assume that there exists a negative curve $f \in [{\bmp_{a,b,c}}^{(r)}]_d$ for some pairwise coprime positive integers $a$, $b$, $c$.
Then we know $\dim_k[{\bmp_{a,b,c}}^{(r)}]_d = 1$ by the same reason as Proposition~\ref{Fact3.3} (2) (a).
Since $[{\bmp_{a,b,c}}^{(r)}]_d$ is defined in $[S_{a,b,c}]_d$ by $\frac{r(r+1)}{2}$ linear equations, we know $\dim_k[S_{a,b,c}]_d \le \frac{r(r+1)}{2} + 1$.
Let $\varphi$ be the $r$-nct corresponding to $f$ as in Proposition~\ref{Prop3.2}.
Let $P$ be the convex hull of $dP_{a,b,c} \cap \bZ^2$.
Then we have
\[
P_\varphi \subset P \subset dP_{a,b,c} .
\]
Usually it is very difficult to determine $P_\varphi$.

Here assume ${}^\#(dP_{a,b,c} \cap \bZ^2) = \dim_k[S_{a,b,c}]_d = \frac{r(r+1)}{2} + 1$.
Then $P$ contains just $\frac{r(r+1)}{2} + 1$ lattice points.
Since $|P| \le |dP_{a,b,c}| < r^2/2$,
we know the number of lattice points in the boundary of $P$ is bigger than or equal to $r+1$ and
that in the interior of $P$ is less than or equal to $r(r-1)/2$ by Pick's theorem.
Since the interior of $P_\varphi$ is contained in  that of $P$,
the condition (9) in Theorem~\ref{Thm3.6} is satisfied for $\varphi$ by Lemma~\ref{Lemma8.9}.
\end{rm}
\end{Remark}

\begin{Remark}\label{rational}
\begin{rm}
Assume that $k$ is algebraically closed and there exists a negative curve $f \in [{\bmp_{a,b,c}}^{(r)}]_d$ for some pairwise coprime positive integers $a$, $b$, $c$.
Let $\varphi$ be the $r$-nct corresponding to $f$ as in Proposition~\ref{Prop3.2}.
Then the negative curve $C$ in $Y_{\Delta_{a,b,c}}$ is birational to $C_\varphi$ in $Y_{\Delta_{\varphi}}$.
If the condition (9) in Theorem~\ref{Thm3.6} is satisfied for $\varphi$,
we know that $C$ is a rational curve.

If $r=1$, it is easy to check $C \simeq \bP^1_k$
(e.g.\ Lemma~3.2 in \cite{KN}).

Suppose $r=2$.
Since the unique $2$-nct satisfies (1) in Theorem~\ref{Thm3.6} (cf.\ Example~\ref{ExampleOfNct} (2),
Example~\ref{ExampleOfGGK}),
$C$ is a rational curve.
It is easy to see that $C$ is singular if and only if ${}^\#(d {P_{a,b,c}}^\circ \cap \bZ^2) > 1$.
There are many examples such that $C \simeq \bP^1_k$ (e.g.\ 
$(a,b,c) = (3,7,8)$, $(16,97,683), \ldots$).
There are also many examples such that $C$ is singular (e.g.\ $(a,b,c) = (5,77,101)$, $(107,159,173), \ldots$).
In the case where $(a,b,c) = (3,7,8)$, $(5,77,101)$, $R'_s(\bmp_{a,b,c})$ is Noetherian.
In the case where $(a,b,c) = (16,97,683)$, $(107,159,173)$, $R'_s(\bmp_{a,b,c})$ is not Noetherian. 

Assume that $r = 3$.
In this case any $3$-nct satisfies (1) in Theorem~\ref{Thm3.6} over any field $k$.
In fact, when ${\rm ch}(k) \neq 2$, we can prove it in the same way as Example~\ref{ExampleOfGGK}.
Assume that ${\rm ch}(k) = 2$.
For $\varphi_3$ in Example~\ref{123-nct}, we can also prove it in the same way as Example~\ref{ExampleOfGGK}.
Consider $\varphi'_3$ in Example~\ref{123-nct}.
Remark that $P_{\varphi'_3}$ is a triangle and the Picard number of $Y_{\Delta_{\varphi'_3}}$ is $2$.
Since $C_{\varphi'_3}^2 < 0$ and $-K_{Y_{\Delta_{\varphi'_3}}}.C_{\varphi'_3} = B_{\varphi'_3} - 3 = 0$, 
we know that $-K_{Y_{\Delta_{\varphi'_3}}}$ is nef and big.

If $r = 4$, then (9) in Theorem~\ref{Thm3.6} is satisfied over any field $k$.
If $r = 4$ and ${\rm ch}(k) = 0$, then (7) in Theorem~\ref{Thm3.6} is satisfied.

If $r = 5$ and ${\rm ch}(k)=0$, then one can prove that (9) in Theorem~\ref{Thm3.6} is satisfied.

The author does not know any example that 
the condition (7) in Theorem~\ref{Thm3.6}  is not satisfied in the case where
${\rm ch}(k) = 0$.

From the above, if $r \le 4$, then $C$ is rational.
If $r = 5$ and ${\rm ch}(k)=0$, then $C$ is rational.

In many case $\dim_k[S_{a,b,c}]_d = \frac{r(r+1)}{2} + 1$ is satisfied.
When this is the case, $C$ is rational since (9) in Theorem~\ref{Thm3.6} is satisfied (cf.\ Remark~\ref{SMC}).
The author knows a few examples that $\dim_k[S_{a,b,c}]_d < \frac{r(r+1)}{2} + 1$.
See the next example.
\end{rm}
\end{Remark}

\begin{Example}\label{(8,15,43)}
\begin{rm}
Suppose that $k$ is of characteristic $0$ and $(a,b,c) = (8,15,43)$.
Using a computer, we know that there exists a negative curve $f \in [{\bmp_{a,b,c}}^{(9)}]_{645}$.
Let $\varphi$ be the corresponding $9$-nct  as in Proposition~\ref{Prop3.2}.
Then $P_\varphi$ is the following pentagon:

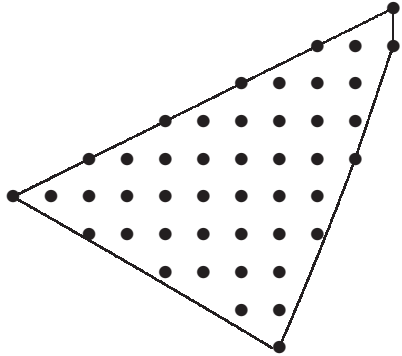
\[
{
\setlength\unitlength{1truecm}
  \begin{picture}(5,4.5)(0,-2)
\qbezier (0.1,0.1) (1.1,0.6) (5.1,2.6)
\qbezier (0.1,0.1) (1.85,-0.9) (3.5,-1.9)
\qbezier (3.6,-1.9) (4.1,-0.75) (4.6,0.6)
\qbezier (4.6,0.6) (4.85,1.35) (5.1,2.1)
\qbezier (5.1,2.1) (5.1,2.4) (5.1,2.6)
\put(0,0){$\bullet$}
\put(0.5,0){$\bullet$}
\put(1,0.5){$\bullet$}
\put(1,0){$\bullet$}
\put(1,-0.5){$\bullet$}
\put(1.5,0.5){$\bullet$}
\put(1.5,0){$\bullet$}
\put(1.5,-0.5){$\bullet$}
\put(2,1){$\bullet$}
\put(2,0.5){$\bullet$}
\put(2,0){$\bullet$}
\put(2,-0.5){$\bullet$}
\put(2,-1){$\bullet$}
\put(2.5,1){$\bullet$}
\put(2.5,0.5){$\bullet$}
\put(2.5,0){$\bullet$}
\put(2.5,-0.5){$\bullet$}
\put(2.5,-1){$\bullet$}
\put(3,1.5){$\bullet$}
\put(3,1){$\bullet$}
\put(3,0.5){$\bullet$}
\put(3,0){$\bullet$}
\put(3,-0.5){$\bullet$}
\put(3,-1){$\bullet$}
\put(3,-1.5){$\bullet$}
\put(3.5,1.5){$\bullet$}
\put(3.5,1){$\bullet$}
\put(3.5,0.5){$\bullet$}
\put(3.5,0){$\bullet$}
\put(3.5,-0.5){$\bullet$}
\put(3.5,-1){$\bullet$}
\put(3.5,-1.5){$\bullet$}
\put(3.5,-2){$\bullet$}
\put(4,2){$\bullet$}
\put(4,1.5){$\bullet$}
\put(4,1){$\bullet$}
\put(4,0.5){$\bullet$}
\put(4,0){$\bullet$}
\put(4,-0.5){$\bullet$}
\put(4.5,2){$\bullet$}
\put(4.5,1.5){$\bullet$}
\put(4.5,1){$\bullet$}
\put(4.5,0.5){$\bullet$}
\put(5,2.5){$\bullet$}
\put(5,2){$\bullet$}
  \end{picture}
}
\]

\noindent
It satisfies $B_\varphi + I_\varphi = \frac{r(r+1)}{2}= 45$,
$B_\varphi = r =9$, $I_\varphi = \frac{r(r-1)}{2}=36$.
In this case, $P_{a,b,c}$ is the triangle with three edges having slopes
$1/2$, $-4/7$ and $5/2$.
Since $a+b+c$ is less than $d/r$ in this case,  $-K_{Y_{\Delta_{a,b,c}}}$ is neither nef nor big.
Since we have a birational surjective map $Y_{\Delta_\varphi} \rightarrow Y_{\Delta_{a,b,c}}$,
we know that $-K_{Y_{\Delta_\varphi}}$ is neither nef nor big.
The author does not know whether ${\rm Cox}(Y_{\Delta_\varphi})$ is Noetherian or not.
Since the condition (7) in Theorem~\ref{Thm3.6} is satisfied,
the extended symbolic Rees ring $R'_s(\bmp_\varphi)$ is Noetherian.
In this case, generators of one-dimensional cones of the fan $\Delta_\varphi$ is 
\[
\left(
\begin{array}{c}
\bma_1 \\ \bma_2 \\ \bma_3 \\ \bma_4 \\ \bma_5
\end{array}
\right)
=
\left(
\begin{array}{cc}
-1 & 0 \\ -3 & 1 \\ -5 & 2 \\ 4 & 7 \\ 1 & -2
\end{array}
\right) .
\]
Therefore the morphism of monoids (\ref{grading}) is 
\[
\bZ^3 \stackrel{
\left(
\begin{array}{ccccc}
1 & -2 & 1 & 0 & 0 \\ 25 & -7 & 0 & 1 & 0 \\ -5 & 2 & 0 & 0 & 1
\end{array}
\right)}{
\longleftarrow} \bZ^5 \supset (\bNo)^5 .
\]
Let $I$ be the kernel of the $k$-algebra homomorphism $k[x_1,x_2,x_3,x_4,x_5] \rightarrow k[y_1^{\pm 1},y_2^{\pm 1},y_3^{\pm 1}]$ 
($x_1 \mapsto y_1y_2^{25}y_3^{-5}$, $x_2 \mapsto y_1^{-2}y_2^{-7}y_3^2$, $x_3 \mapsto y_1$, $x_4 \mapsto y_2$, $x_5 \mapsto y_3$).
Then we have ${\rm Cox}(Y_{\Delta_\varphi})=R'_s(I)$ as in Remark~\ref{coxY}.

Using a computer, we know that there exists a negative curve $f \in [P^{(18)}]_{1617}$
in the case where $(a,b,c) = (5,33,49)$ over a field $k$ of characteristic $0$.
Let $\varphi$ be the corresponding $18$-nct as in Proposition~\ref{Prop3.2}.
Let $P$ be the convex hull of $1617P_{5,33,49} \cap \bZ^2$ as in Remark~\ref{SMC}.
The number of lattice points in the boundary of $P$ is $r=18$ and
that in the interior of $P$ is $r(r-1)/2=18\times17/2$.
Since $P_\varphi$ is contained in $P$,  the number of lattice points in the interior of $P_\varphi$ is $r(r-1)/2=18\times17/2$ by Lemma~\ref{Lemma8.9}.
Hence the condition (9) in Theorem~\ref{Thm3.6} is satisfied for $\varphi$.
In this case, $-K_{Y_{\Delta_\varphi}}$ is neither nef nor big.
The author does not know whether ${\rm Cox}(Y_{\Delta_\varphi})$ is Noetherian or not.
\end{rm}
\end{Example}

\begin{Remark}\label{K+C}
\begin{rm}
Let $\varphi$ be an $r$-nct.
Then we have 
$h^0(Y_{\Delta_\varphi}, {\cal O}_{Y_{\Delta_\varphi}}(-K_{Y_{\Delta_\varphi}} + C_\varphi)) = p_a(C_\varphi)$,
$h^1(Y_{\Delta_\varphi}, {\cal O}_{Y_{\Delta_\varphi}}(-K_{Y_{\Delta_\varphi}} + C_\varphi)) = 0$
and $I_\varphi = \frac{r(r-1)}{2} + p_a(C_\varphi)$.

Let $f \in [{\bmp_{a,b,c}}^{(r)}]_d$ be the negative curve.
Then we have 
$h^0(Y_{\Delta_{a,b,c}}, {\cal O}_{Y_{\Delta_{a,b,c}}}(-K_{Y_{\Delta_{a,b,c}}} + C)) = p_a(C)$,
$h^1(Y_{\Delta_{a,b,c}}, {\cal O}_{Y_{\Delta_{a,b,c}}}(-K_{Y_{\Delta_{a,b,c}}} + C)) = 0$
and the number of lattice points in the interior of $dP_{a,b,c}$ is $\frac{r(r-1)}{2} + p_a(C)$.
\end{rm}
\end{Remark}

\vspace{3mm}

\noindent
\begin{tabular}{l}
Kazuhiko Kurano \\
Department of Mathematics \\
Faculty of Science and Technology \\
Meiji University \\
Higashimita 1-1-1, Tama-ku \\
Kawasaki 214-8571, Japan \\
{\tt kurano@meiji.ac.jp} \\
{\tt http://www.isc.ac.jp/\~{}kurano}
\end{tabular}

\end{document}